\documentclass[a4paper,dvipsnames]{article}

\usepackage{geometry}
\usepackage[british]{babel}
\usepackage{amsmath}
\usepackage{amssymb}
\usepackage{amsthm}
\usepackage{upgreek}
\usepackage{bm}
\usepackage{cite}

\usepackage[colorlinks=true,allcolors=purple]{hyperref}
\usepackage{cleveref}
	\crefformat{equation}{#2(#1)#3}
\usepackage{dsfont}
\usepackage{enumitem}
\usepackage[OT2,OT1]{fontenc}
\usepackage[cal=cm,scr=boondoxo]{mathalfa}
\usepackage{pifont}
\usepackage{stmaryrd}
\usepackage{textcomp}
\usepackage{tikz}
	\usetikzlibrary{arrows}
	\usetikzlibrary{patterns}
\usepackage{tikz-cd}
\usepackage[textwidth=3cm,colorinlistoftodos]{todonotes}
\usepackage{xfrac}
\usepackage[UKenglish]{isodate}

\numberwithin{equation}{section}			
\swapnumbers						

\newcommand\cyr{
\renewcommand\rmdefault{wncyr}%
\renewcommand\sfdefault{wncyss}%
\renewcommand\encodingdefault{OT2}%
\normalfont
\selectfont}
\DeclareTextFontCommand{\textcyr}{\cyr}

\newcounter{Enum}				
\newenvironment{Enumerate}{\begin{enumerate}[label={\rm({\roman*})}]}{\end{enumerate}}
\newcommand{\Enumref}[1]{{\setcounter{Enum}{#1}{\rm(\roman{Enum})}}}

\newcommand{\descriptionlabelsave}{}		
\newenvironment{Itemize}{%
	\renewcommand{\descriptionlabelsave}{\descriptionlabel}\renewcommand{\descriptionlabel}{$\triangleright$}%
	\begin{description}[leftmargin=15pt,itemindent=-5.2pt]}{%
	\end{description}\renewcommand{\descriptionlabel}{\descriptionlabelsave}}

\newcounter{StepsCount}				
\newenvironment{Elist}{%
	\begin{list}{\ding{\value{StepsCount}}}{\usecounter{StepsCount} \leftmargin=0pt \labelwidth=12pt \itemindent=\labelwidth%
	\itemsep=5pt\listparindent=\parindent} \setcounter{StepsCount}{191}}{\end{list}}
\newcounter{StepsRefCount}

\newenvironment{Ilist}{
	\begin{list}{$\triangleright$}{\leftmargin=0pt \labelwidth=11pt \itemindent=\labelwidth%
	\itemsep=5pt\listparindent=\parindent}}{\end{list}}

\theoremstyle{plain}
	\newtheorem{lemma}{Lemma}[section]
	\newtheorem{proposition}[lemma]{Proposition}
	\newtheorem{theorem}[lemma]{Theorem}
	\newtheorem{corollary}[lemma]{Corollary}
	\newcommand{\GenericTheoremName}{}\newtheorem{generictheorem}[lemma]{\GenericTheoremName}
\theoremstyle{definition}
	\newtheorem{definition}[lemma]{Definition}
	\newcommand{\GenericDefinitionName}{}\newtheorem{genericdefinition}[lemma]{\GenericDefinitionName}
\theoremstyle{remark}
	\newtheorem{remark}[lemma]{Remark}
	\newtheorem{example}[lemma]{Example}
	\newcommand{\GenericRemarkName}{}\newtheorem{genericremark}[lemma]{\GenericRemarkName}
\newenvironment{Proposition}{\begin{proposition}}{\par\noindent\centerline{\rule{5em}{1pt}}\end{proposition}}
\newenvironment{Theorem}{\begin{theorem}}{\par\noindent\centerline{\rule{5em}{1pt}}\end{theorem}}
\newenvironment{GenericTheorem}[1]
	{\renewcommand{\GenericTheoremName}{#1}\begin{generictheorem}}{\par\noindent\centerline{\rule{5em}{1pt}}\end{generictheorem}}
\newenvironment{Definition}{\begin{definition}}{\par\noindent\centerline{\rule{5em}{1pt}}\end{definition}}
\newenvironment{Remark}{\begin{remark}}{\par\noindent\centerline{\rule{5em}{1pt}}\end{remark}}
\newenvironment{Example}{\begin{example}}{\par\noindent\centerline{\rule{5em}{1pt}}\end{example}}

\newcommand{\mc}[1]{{\mathcal{#1}}}			
\newcommand{\ms}[1]{{\mathscr{#1}}}			
\newcommand{\mf}[1]{{\mathfrak{#1}}}			
\newcommand{\bb}[1]{{\mathbb{#1}}}			
\newcommand{\ov}{\overline}				
\newcommand{\mr}{\mathring}				
\newcommand{\wt}{\widetilde}				

\DeclareMathOperator{\RE}{Re}				
\renewcommand{\Re}{\RE}
\DeclareMathOperator{\IM}{Im}				
\renewcommand{\Im}{\IM}

\DeclareMathOperator{\Arccot}{Arccot}			

\newcommand{\smmatrix}[4]{\Bigl(			
\begin{smallmatrix}
\hspace*{-0.2ex} #1 \hspace*{0.2ex} & \hspace*{0.2ex} #2 \hspace*{-0.2ex}
\\[0.5ex]
\hspace*{-0.2ex} #3 \hspace*{0.2ex} & \hspace*{0.2ex} #4 \hspace*{-0.2ex}
\end{smallmatrix}
\Bigr)}

\newcommand{\Dummy}{\text{\textvisiblespace\kern1pt}}	
\newcommand{\Smallo}{{\rm o}}				

\newcommand{\DS}{\colon\mkern3mu}			
\newcommand{\DP}{{.\kern7pt}}				
\newcommand{\DF}{\colon}				
\newcommand{\DE}{\mathrel{\mathop:}=}			
\newcommand{\ED}{=\mathrel{\mathop:}}			
\newcommand{\DD}{\mkern4mu\mathrm{d}}			

\DeclareMathOperator{\LP}{LP}				
\DeclareMathOperator{\Loc}{loc}				
\DeclareMathOperator{\sgn}{sgn}				
\newcommand{\dd}{\hat{d}}				

\newcommand{\Ham}{{\bb H}}				
\newcommand{\HamInt}[1]{
	\Ham_{#1}}
\newcommand{\NHam}{{\bb H^{1}}}				
\newcommand{\NHamInt}[1]{
	\bb H^{1}_{#1}}
\newcommand{\scHam}{{\bb H}^{\textsf{cs}}}			
\newcommand{\scHamInt}[1]{
	\scHam_{#1}}

\newcommand{\genResc}[2]{\ms A_{#1}{#2}}			
\newcommand{\trResc}[1]{\tau_{#1}}			

\DeclareMathOperator{\tr}{tr}				
\DeclareMathOperator{\Dist}{dist}			

\newcommand{\hatM}{\rule{0ex}{1ex}\hspace{0.6ex}\widehat{\rule[1.5ex]{1.5ex}{0ex}}\hspace{-2.2ex}M}
\newcommand{\hatW}{\rule{0ex}{1ex}\hspace{0.6ex}\widehat{\rule[1.5ex]{1.5ex}{0ex}}\hspace{-2.0ex}W}

\newcounter{counter_a}
\newenvironment{enumdash}{\begin{list}{{\rm(\roman{counter_a}$'$)}}%
{\usecounter{counter_a}
\setlength{\itemsep}{0.7ex}\setlength{\topsep}{0.7ex}
\setlength{\leftmargin}{6ex}\setlength{\labelwidth}{6ex}}}{\end{list}}

\newcounter{counter_b}
\newenvironment{enumtripledash}{\begin{list}{{\rm(\roman{counter_b}$'''$)}}%
{\usecounter{counter_b}
\setlength{\itemsep}{0.7ex}\setlength{\topsep}{0.7ex}
\setlength{\leftmargin}{7ex}\setlength{\labelwidth}{7ex}}}{\end{list}}

\begin{document}

\begin{flushleft}
	{\Large\textbf{Canonical systems whose Weyl coefficients have \\[1ex] dominating real part}}
	\\[5mm]
	\textsc{
	Matthias Langer
	\,\ $\ast$\,\
	Raphael Pruckner
	\,\ $\ast$\,\
	Harald Woracek
		\hspace*{-14pt}
		\renewcommand{\thefootnote}{\fnsymbol{footnote}}
		\setcounter{footnote}{2}
		\footnote{The second and third authors were supported by the project P~30715-N35
			of the Austrian Science Fund (FWF).  The third author was supported by the
			joint project I~4600 of the Austrian Science Fund (FWF) and the
			Russian foundation of basic research (RFBR).}
		\renewcommand{\thefootnote}{\arabic{footnote}}
		\setcounter{footnote}{0}
	}
\end{flushleft}
	\vspace*{3ex}
	{\small
	\textbf{Abstract.}
		For a two-dimensional canonical system $y'(t)=zJH(t)y(t)$ on the half-line $(0,\infty)$
		whose Hamiltonian $H$ is a.e.\ positive semi-definite, denote by $q_H$ its Weyl coefficient.
		De~Branges' inverse spectral theorem states that the assignment
		$H\mapsto q_H$ is a bijection between Hamiltonians (suitably normalised)
		and Nevanlinna functions.

		The main result of the paper is a criterion when the singular integral
		of the spectral measure, i.e.\ $\Re q_H(iy)$,
		dominates its Poisson integral $\Im q_H(iy)$ for $y\to+\infty$.
		Two equivalent conditions characterising this situation are provided.
		The first one is analytic in nature, very simple, and explicit in terms of the
		primitive $M$ of $H$.
		It merely depends on the relative size of the off-diagonal entries of $M$ compared
		with the diagonal entries.
		The second condition is of geometric nature and technically more complicated.
		It involves the relative size of the off-diagonal entries of $H$,
		a measurement for oscillations of the diagonal
		of $H$, and a condition on the speed and smoothness of the rotation of $H$.
	\\[3mm]
	}

\vspace*{2ex}

%


%
%
%
\section{Introduction}

We investigate the spectral theory of two-dimensional \emph{canonical systems}
\begin{equation}\label{G60}
	y'(t)=zJH(t)y(t),\qquad t\in(a,b),
\end{equation}
where $-\infty<a<b\leq\infty$, $z\in\bb C$ is the spectral parameter, $J$ is the
symplectic matrix $J\DE\smmatrix 0{-1}10$, and $H$ is the \emph{Hamiltonian} of the system.
We deal with systems whose Hamiltonian satisfies
\begin{Ilist}
\item $H(t)\in\bb R^{2\times 2}$ and $H(t)\geq 0$ a.e.;
\item for all $c\in(a,b)$ we have $\int_a^c\tr H(s)\DD s<\infty$;
\item $H(t)\ne 0$ a.e.
\end{Ilist}
We further assume that $H$ is in the limit point case at the
right endpoint $b$, i.e.\
\begin{equation}\label{G1}
	\int_a^b\tr H(s)\DD s=\infty.
\end{equation}
A central role in the theory of such equations is played by the Weyl coefficient $q_H$
associated with $H$.
For Sturm--Liouville equations its construction goes back to H.~Weyl \cite{weyl:1910}.
Let us recall the definition of $q_H$ for canonical systems.
To this end, let $W(t,z)$ be the (transpose of) the fundamental solution of the
system \eqref{G60},
i.e.\ the unique $2\times2$-matrix-valued solution of the initial value problem
\[
	\left\{
	\begin{array}{l}
		\dfrac{\partial}{\partial t}W(t,z)J = zW(t,z)H(t),\qquad t\in[a,b),
		\\[2.5ex]
		W(a,z)=I.
	\end{array}
	\right.
\]
Note that the transposes of the rows of $W$ are solutions of \eqref{G60},
and let us write $W(t,z)=\smmatrix{w_{11}(t,z)}{w_{12}(t,z)}{w_{21}(t,z)}{w_{22}(t,z)}$.
If \cref{G1} is satisfied, then the following limit exists and
is independent of $\zeta$ in the closed upper half-plane $\bb C^+\cup\bb R$:
\[
	q_H(z) \DE \lim_{t\to b}\frac{w_{11}(t,z)\zeta+w_{12}(t,z)}{w_{21}(t,z)\zeta+w_{22}(t,z)}\,,
	\qquad z\in\bb C\setminus\bb R;
\]
the function $q_H$ is called \emph{Weyl coefficient} associated with the Hamiltonian $H$.
It is a Nevanlinna function or identically equal to $\infty$
(when $h_2(t)=0$ for a.e.\ $t\in(a,b)$);
a \emph{Nevanlinna function}%
\footnote{%
	Sometimes in the literature the terminology \emph{Herglotz function} is used instead.
}
is a function that is analytic in $\bb C\setminus\bb R$ and
satisfies $q_H(\ov z)=\ov{q_H(z)}$ and $\Im q_H(z)\cdot\Im z\ge 0$ for all $z$.
The significance of the Weyl coefficient is that the measure $\mu$ in its
Herglotz integral representation
\[
	q_H(z) = \alpha+\beta z+\int_{\bb R}\Bigl(\frac 1{t-z}-\frac t{1+t^2}\Bigr)\DD\mu
\]
is a spectral measure for the differential operator constructed from the
equation \cref{G60}
(when $\beta>0$, this differential operator is actually multi-valued and one can
include a point mass at infinity with mass $\beta$).

A famous theorem by L.~de~Branges \cite{debranges:1968} says that the
assignment $H\mapsto q_H$ establishes a bijective correspondence between the set of
all suitably normalised Hamiltonians on the one hand, and the set of all Nevanlinna functions
on the other hand.  In view of de~Branges' correspondence, it is a natural task to
translate properties from $H$ to $q_H$ (i.e.\ \emph{direct spectral relations})
and vice versa from $q_H$ to $H$ (i.e.\ \emph{inverse spectral relations}).
In the best case one can go both ways.
For illustration, let us mention two examples of such theorems.
It is possible to explicitly characterise those Hamiltonians $H$ for which $q_H$
has an analytic continuation to $\bb C\setminus[0,\infty)$, see \cite{winkler:1998},
or those Hamiltonians for which $q_H$ has a meromorphic continuation to all of $\bb C$,
see \cite{romanov.woracek:ideal}.
The first result characterises that the differential operator associated with \cref{G60}
is non-negative, the second one that it has discrete spectrum.

In the present paper we prove a direct and inverse spectral relation of a different kind.
It belongs to a family of results which relate the behaviour of $H$ locally at the
left endpoint $a$ with the behaviour of $q_H$ when $z$ tends to $+i\infty$;
for physical reasons one also speaks of the \emph{high-energy behaviour} of $q_H$.
Recall that the behaviour of $\Im q_H(iy)$ at $+\infty$ is related to the behaviour of the
spectral measure at $\pm\infty$; see, e.g.\ \cite[Section~4]{langer.pruckner.woracek:heniest-arXiv}.
Our main result is \Cref{G14} stated further below, where we characterise those
Hamiltonians $H$ for which%
\footnote{%
	We use the notation ``$f\ll g$'' for $f/g\to 0$.
}%
\begin{equation}\label{G61}
	\Im q_H(iy)\ll |q_H(iy)|,\qquad y\to+\infty,
\end{equation}
i.e.\ those Hamiltonians for which the singular integral $\Re q_H(z)$ of the
spectral measure strictly dominates the Poisson integral $\Im q_H(z)$.

In our theorem, where \cref{G61} is listed as item \Enumref{1}, we give two different
conditions on $H$, called \Enumref{2} and \Enumref{3}, which are both equivalent to \cref{G61}.
Condition \Enumref{2} is analytic in nature, very simple, and explicit in terms of
the primitive
\[
	M(t)\DE\int_a^t H(s)\DD s
\]
of $H$, which is a nonnegative and nondecreasing matrix function.
It says that, locally at $a$, the off-diagonal entries of $M(t)$ should be as large as its diagonal entries.
Condition \Enumref{3} is of geometric nature and somewhat more complicated.
It involves the relative size of the off-diagonal entries of $H$
compared with the diagonal entries, a measurement for oscillations of the diagonal
of $H$, and a condition on the speed and smoothness of the ``rotation'' of $H$.

From a function-theoretic perspective, the behaviour exhibited by \cref{G61} is
rather peculiar.  For every Nevanlinna function $q$ one has that for
(in a measure-theoretic sense) most points on the boundary of the open upper half-plane
(including $+i\infty$) condition \cref{G61} fails; see \cite{poltoratski:2003} and
recall that real and imaginary parts are comparable on approaching almost every point of
the absolutely continuous spectrum.
On the other hand, for a certain subclass of Nevanlinna functions it holds that for
(in a topological sense) many boundary points \cref{G61} holds,
cf.\ \cite[Theorem~1]{donaire:2001} where one uses a curve that
approaches the boundary tangentially.
Neither of these statements has any implication for a single boundary point
(in our case $+i\infty$).
The condition \Enumref{3} in \Cref{G14} is a very strong restriction on $H$.
Hence, one message of \Cref{G14} is that \cref{G61}, i.e.\ strict dominance of the
singular integral at a specific boundary point, is a rather rare phenomenon.

Our interest in the class of Hamiltonians with \cref{G61} originates from the recent result
\cite[Theorem~1.1]{langer.pruckner.woracek:heniest-arXiv}.  In this theorem we showed that,
for every Hamiltonian $H$, the following estimates%
\footnote{%
	We write ``$f\lesssim g$'' for $\exists c>0\DP f\leq cg$,
	and ``$f\asymp g$'' for $f\lesssim g\wedge g\lesssim f$.
}%
\begin{equation}\label{G7}
	|q_H(iy)|\asymp A_H(y) \qquad\text{and}\qquad
	L_H(y)\lesssim\Im q_H(iy)\lesssim A_H(y) \qquad\text{for $y\ge 1$},
\end{equation}
hold, where $L_H(y)$ and $A_H(y)$ are certain functions defined explicitly in terms
of the primitive $M(t)$, and the constants in ``$\asymp$'' and ``$\lesssim$'' are
independent of $H$; we recall details in \Cref{G84}.
The question arises whether the lower bound $L_H(y)$ is sharp.
The equivalence of \cref{G61} with \Cref{G14}\,\Enumref{2} says that on a qualitative level
the answer is affirmative: we have
\[
	\Im q_H(iy)\ll|q_H(iy)|\quad\Leftrightarrow\quad L_H(iy)\ll A_H(iy).
\]
It is an open problem if there is a quantitative relation between
$\Im q_H(iy)$ and $L_H(iy)$ (assuming that $\Im q_H(iy)\ll|q_H(iy)|$ and
thinking up to universal multiplicative constants).
This seems to be a rather involved question, and we expect that the equivalence of \cref{G61}
with \Cref{G14}\,\Enumref{3} will be of help to attack it.

Let us give a brief overview of the contents of the paper.
In the remainder of the Introduction we formulate the main theorem,
\Cref{G14}, and a sequence variant, \Cref{G15}, and provide an illustrative example.
In \Cref{G56} we provide some preliminaries and set up notation.
\Cref{G123} contains the proof of the equivalence of (i) and (ii) in our main results.
\Cref{G55} contains preparations for the proof of the equivalence with (iii),
which is then carried out in \Cref{G83}.
Finally, in \Cref{G98} we consider the situation when the diagonal entries
of $H$, or their primitives, are regularly varying.

\subsection*{Formulation of the main theorem}

We formulate our main theorem for Hamiltonians that satisfy
\begin{Ilist}
\item
	$a=0$, $b=\infty$;
\item
	neither of the diagonal entries of $H$ vanishes a.e.\ on some interval starting
	at the left endpoint~$0$.
\end{Ilist}
Both assumptions are no loss in generality, and are only imposed for simplicity.
The first one can always be achieved by a change of the independent variable in
equation \cref{G60}, and changes of variable do not alter the Weyl coefficient;
see \Cref{G70}.
The second condition excludes some exceptional cases where there is nothing to investigate:
if it is not satisfied, then $\lim_{y\to\infty}\frac{\Im q_H(iy)}{|q_H(iy)|}=1$;
we provide more details in \Cref{G70,G71}.

Throughout the paper we write
\begin{equation}\label{G54}
	H(t)=
	\begin{pmatrix}
		h_1(t) & h_3(t)
		\\[0.5ex]
		h_3(t) & h_2(t)
	\end{pmatrix}
	,\qquad m_j(t)\DE\int_0^t h_j(s)\DD s,\quad j=1,2,3;
\end{equation}
sometimes we write $M(H,t)$ and $m_i(H,t)$ instead of $M(t)$ and $m_i(t)$ respectively
to indicate the dependence on $H$.
Moreover, $\lambda$ denotes the Lebesgue measure.

Next, we have to introduce some notation which looks a bit technical on first sight,
but actually is not.  The intuition behind these quantities is discussed in \Cref{G63} below.
The functions are well defined because $h_3(t)^2\le h_1(t)h_2(t)$ for a.e.\ $t>0$
and $m_1(t),m_2(t)>0$ for all $t>0$; the latter follows from the assumption
that neither of the diagonal entries of $H$ vanishes a.e.\ on an interval starting at $0$.
Set
\begin{align}
	\sigma_H(t)\DE &\,
	\begin{cases}
		\dfrac{|h_3(t)|}{\sqrt{h_1(t)h_2(t)}\,} &\text{if}\ h_3(t)\ne 0,
		\\[3ex]
		0 &\text{otherwise},
	\end{cases}
	\label{G57}
	\\[1ex]
	\pi_{H,s}(t)\DE &\,
	\begin{cases}
		\displaystyle\raisebox{5pt}{${\sgn\bigl(h_3(st)\bigr)}\dfrac{h_2(st)}{h_1(st)}$}
		\Bigg/
		\raisebox{-5pt}{$\dfrac{m_2(s)}{m_1(s)}$} &\text{if}\ h_3(st)\ne 0,
		\\[3ex]
		0 &\text{otherwise},
	\end{cases}
	\label{G86}
	\\[1ex]
	\mf t_s(t)\DE &\,
	\frac{m_1(st)}{m_1(s)}+\frac{m_2(st)}{m_2(s)},
	\label{G76}
\end{align}
where $s>0$ is a parameter.

Note that, for each fixed $s>0$, the function $\mf t_s$ is absolutely continuous
and its derivative
\[
	\mf t_s'(t)=\frac s{m_1(s)}h_1(st)+\frac s{m_2(s)}h_2(st)
\]
is positive a.e.  Furthermore, 
$\mf t_s(0)=0$
and $\lim_{t\to\infty}\mf t_s(t)=\infty$;
the latter follows from the relation $m_1(st)+m_2(st)=\int_0^{st}\tr H(x)\DD x\to\infty$
as $t\to\infty$ by assumption.
Thus $\mf t_s$ is an increasing bijection from $[0,\infty)$ onto itself with
absolutely continuous inverse function.

Now we are in position to state our main theorem.

\begin{Theorem}\label{G14}
	Let $H$ be a Hamiltonian defined on the interval $(0,\infty)$ such that \textup{\cref{G1}} holds
	and neither $h_1$ nor $h_2$ vanishes a.e.\ on some neighbourhood of the left endpoint $0$.
	Then the following statements are equivalent.
	\begin{Enumerate}
	\item
		Relation \textup{\cref{G61}} holds, i.e.\
		\begin{equation}\label{G92}
			\lim_{y\to\infty}\frac{\Im q_H(iy)}{|q_H(iy)|}=0.
		\end{equation}
	\item
		We have
		\begin{equation}\label{G64}
			\lim_{t\to 0}\frac{\det M(t)}{m_1(t)m_2(t)}=0.
		\end{equation}
	\item
		For all $T\in(0,\infty)$, all $\gamma\in[0,1)$, and all open intervals
		$I,J\subseteq\bb R\setminus\{0\}$ with $\ov I\cap\ov J=\emptyset$ and
		at least one of $I$ and $J$ being bounded, the following limit relations hold:
		\begin{align}
			& \lim_{s\to 0}\Bigl[\lambda\Bigl(
			(0,T)\cap\mf t_s\bigl(\tfrac 1s\sigma_H^{-1}([0,\gamma])\bigr)
			\Bigr)\Bigr]=0,
			\label{G65}
			\\[1ex]
			& \lim_{s\to 0}\Bigl[
			\lambda\Bigl((0,T)\cap\mf t_s(\pi_{H,s}^{-1}(I))\Bigr)
			\cdot\lambda\Bigl((0,T)\cap\mf t_s(\pi_{H,s}^{-1}(J))\Bigr)\Bigr]=0.
			\label{G16}
		\end{align}
	\end{Enumerate}
\end{Theorem}

\noindent
Under a certain additional assumption, the conditions in \Enumref{3} greatly simplify.
This assumption is quite strong, and will, in many interesting cases, not be satisfied.
Still, in order to understand the nature of \cref{G65} and \cref{G16} and the
proof of their equivalence to \cref{G61}, it is worth stating the following addition.

\begin{GenericTheorem}{Addition to \Cref{G14}}
	Assume that, in addition to the assumptions of \Cref{G14}, the following
	conditions hold:
	\begin{align}
		& \tr H(t) = 1 \qquad  \text{for a.e.} \ t\in(0,\infty),
		\label{G99}
		\\[1ex]
		& \liminf\limits_{t\to 0}\Bigl(\frac{m_1(t)}{t}\cdot\frac{m_2(t)}{t}\Bigr)>0.
		\label{G100}
	\end{align}
	Then the equivalent properties \Enumref{1}, \Enumref{2}, \Enumref{3} in \Cref{G14}
	are further equivalent to the following condition.
	\begin{Enumerate}
	\setcounter{enumi}{3}
	\item
		For all $\gamma$ and $I,J$ as in \Cref{G14}\,\Enumref{3} we have
		\begin{align}
			& \lim_{t\to 0}\biggl[\frac 1t\lambda\Bigl(
			(0,t)\cap\sigma_H^{-1}\bigl([0,\gamma]\bigr)\Bigr)\biggr]=0,
			\label{G62}
			\\[1ex]
			& \lim_{t\to 0}\biggl[
			\frac 1t\lambda\Bigl(
			(0,t)\cap\pi_H^{-1}(I)\Bigr)
			\cdot\frac 1t\lambda\Bigl((0,t)\cap\pi_H^{-1}(J)\Bigr)\biggr]=0,
			\label{G66}
		\end{align}
		where
		\begin{equation}\label{G87}
			\pi_H(t)\DE
			\begin{cases}
				{\sgn\bigl(h_3(t)\bigr)}\dfrac{h_2(t)}{h_1(t)} &\text{if}\ h_3(t)\neq 0,
				\\[2ex]
				0 &\text{otherwise}.
			\end{cases}
		\end{equation}
	\end{Enumerate}
\end{GenericTheorem}

\noindent
Note that \cref{G99} implies that $m_1(t)+m_2(t)=t$.
Hence, by \cite[Theorem~1.1]{langer.pruckner.woracek:heniest-arXiv} (see also \Cref{G85})
we have
\begin{equation}
\label{G24}
\begin{alignedat}{2}
	\liminf\limits_{t\to 0}\Bigl(\frac 1t m_1(t)\cdot\frac 1t m_2(t)\Bigr) > 0
	\qquad&\Leftrightarrow\qquad
	\frac{m_1(t)}{m_2(t)} \asymp 1, \quad && t\to0
	\\[1ex]
	&\Leftrightarrow\qquad
	|q_H(ir)|\asymp 1, \quad && r\to\infty.
\end{alignedat}
\end{equation}
We come to the promised explanation of the conditions in \Enumref{3} (and \Enumref{4}).

\begin{Remark}\label{G63}
	Let us first discuss the simpler conditions \cref{G62} and \cref{G66}.

	The role of $\sigma_H$ is to quantify the relative size of the off-diagonal entries
	of $H$ compared with the diagonal entries.
	Condition \cref{G62} can be written as
	$\lim_{t\to0}[\frac1t\lambda(\{x\in(0,t):1-\sigma_H(x)^2\ge1-\gamma^2\})]=0$,
	or, by rescaling, as
	\[
		\lim_{t\to0}\lambda\bigl(\bigl\{x\in(0,1):1-\sigma_H(tx)^2\ge1-\gamma^2\bigr\}\bigr)=0.
	\]
	The validity of this relation for all $\gamma\in[0,1)$ just says that
	the functions $x\mapsto 1-\sigma_H(tx)^2$ converge to $0$ in measure as $t\to0$.
	Since they are non-negative and bounded by $1$, this is also equivalent
	to the fact that their integrals converge to $0$.
	Note that
	\[
		1-\sigma_H(x)^2
		=
		\begin{cases}
			\dfrac{\det H(x)}{h_1(x)h_2(x)} & \text{if} \ h_3(x)\ne0,
			\\[2ex]
			1 & \text{otherwise}.
		\end{cases}
	\]
	Hence the validity of \cref{G62} for all $\gamma\in[0,1)$ is
	(again by rescaling) equivalent to
	\[
		\lim_{t\to 0}\frac 1t\int_0^t\frac{\det H(x)}{h_1(x)h_2(x)}\DD x=0,
	\]
	where the integrand is understood as equal to $1$ at points where its
	denominator vanishes; this means that the Hamiltonian should be almost
	of zero determinant in the vicinity of the left endpoint $0$ in a
	measure-theoretic sense.

	The role of $\pi_H$ is not so obvious.  It is related to what one may think of as
	``rotation'' of $H$.  To see this, write $H$ in the form
	\begin{equation}\label{G34}
		H(t)=
		\begin{pmatrix}
			1 & \sigma_H(t)
			\\
			\sigma_H(t) & 1
		\end{pmatrix}
		\odot
		\biggl[\binom{\cos\varphi_H(t)}{\sin\varphi_H(t)}
		\binom{\cos\varphi_H(t)}{\sin\varphi_H(t)}^*\,\biggr],
	\end{equation}
	where $\odot$ denotes the Hadamard, i.e.\ entry-wise, product of the $2\times 2$-matrices.
	The first factor takes the relative size of the off-diagonal entries into account;
	the second factor has zero determinant and corresponds to some kind of rotation.
	The factorisation in \cref{G34} is possible, for instance, with
	\begin{equation}\label{G58}
		\varphi_H(t)\DE
		\begin{cases}
			\Arccot\sqrt\frac{h_1(t)}{h_2(t)} &\text{if}\ h_2(t) \ne 0,h_3(t) \ge 0,
			\\[1ex]
			\pi-\Arccot\sqrt\frac{h_1(t)}{h_2(t)} &\text{if}\ h_2(t) \ne 0,h_3(t)<0,
			\\[1ex]
			0 &\text{if}\ h_2(t)=0,
		\end{cases}
	\end{equation}
	where $\Arccot$ is the branch with values in $(0,\pi)$.
	Then
	\begin{equation}\label{G94}
		\pi_H(t) = \sgn\bigl(\tfrac\pi2-\varphi_H(t)\bigr)\cdot\tan^2\varphi_H(t).
	\end{equation}
	Now we map $\varphi_H(t)\in[0,\pi)$ onto the unit circle $\bb T$ by
	setting
	\begin{equation}\label{G59}
		\zeta_H(t) \DE e^{2i\varphi_H(t)}.
	\end{equation}
	We may say\,---\,descriptively\,---\,that $\zeta_H$ is the rotation of $H$.

	The statement \cref{G66} is equivalent to the following statement (see \Cref{G83}):
	there are no two separated arcs on the unit circle, such that, in the vicinity
	of the left endpoint 0, $\zeta_H(t)$ often belongs to one arc and also often
	belongs to the other arc.
	In other words, the Hamiltonian should rotate so slowly that,
	on every interval close to $0$, it looks\,---\,from a measure-theoretic
	viewpoint\,---\,as if its direction were constant;
	see also \Cref{G97}

	The more complicated conditions \cref{G65} and \cref{G16} are weighted and
	rescaled variants of \cref{G62} and \cref{G66}; see \Cref{G73}.
	The role of the function $\mf t_s$ is to take care of heavy oscillations,
	and the purpose of the weight $\frac{m_2(s)}{m_1(s)}$ in the definition of
	$\pi_{H,s}$ is to level out the contributions of the two diagonal entries.
	Moreover, zooming into the vicinity of the left endpoint $0$ is now achieved
	by sending the rescaling parameter $s$ to $0$.

	Let us note that also the relation \cref{G65} can be rewritten in integral form,
	namely as
	\[
		\lim_{s\to 0}\int_0^T\frac{\det H(s\mf t_s^{-1}(t))}{(h_1h_2)(s\mf t_s^{-1}(t)))}
		\DD t=0.
	\]
\end{Remark}

\noindent
To prove \Cref{G14} we show the implications
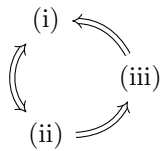
\begin{center}
\begin{tikzcd}[column sep=small, row sep=tiny]
	\textup{\text{(i)}} \arrow[Leftrightarrow, bend right=40]{dd} & &
	\\
	& \textup{\text{(iii)}} \arrow[Rightarrow, bend right=30]{lu}
	\\
	\textup{\text{(ii)}} \arrow[Rightarrow, bend right=30]{ru} & &
\end{tikzcd}
\end{center}
Interestingly, very different methods enter in the proofs of the various implications.
\begin{Itemize}
\item
	The implication ``\Enumref{1}$\Rightarrow$\Enumref{2}'' is a direct consequence of
	\cite[Theorem~1.1]{langer.pruckner.woracek:heniest-arXiv} in the
	form of \Cref{G85} below.
	We recall that this theorem is proved by directly studying Weyl discs and
	estimating the power series coefficients of the fundamental solution of the
	canonical system.
\item
	The proof of ``\Enumref{2}$\Rightarrow$(iii)'' requires an elementary but
	elaborate analysis of the connection between $H$ and its primitive $M$.
	In particular, estimates are proved where the constants are independent
	of the Hamiltonian.
	This is done in \Cref{G55}; see \Cref{G19,G23}.
\item
	To show ``(iii)$\Rightarrow$\Enumref{1}'' and ``(ii)$\Rightarrow$(i)'' we use cluster sets and compactness arguments
	for Hamiltonians endowed with the inverse limit topology of weak $L^1$-topologies
	on finite intervals; see \Cref{G72}.
	Another necessary tool is provided in \Cref{G74}, and a crucial role is taken by
	a weighted variant of Y.~Kasahara's rescaling trick \cite{kasahara:1975},
	which relates the behaviour of $q_H$ towards $i\infty$ with
	weighted rescalings of $H$; see \Cref{G73}.
\end{Itemize}
The proof of ``(ii)$\Rightarrow$(i)'' was included in order to decouple the equivalences between
(i) and (ii), and between (i) and (iii), respectively. This enables reading the proof of ``(i)$\Leftrightarrow$(ii)''
without having to go into the technical details of \Cref{G55}.
We thank a referee for suggesting an argument which makes this possible.

\subsection*{A sequence variant of the theorem}

We can also give a variant of \Cref{G14} where limits are replaced by limits inferior.
It reads as follows.

\begin{Theorem}\label{G15}
	Let $H$ be a Hamiltonian defined on the interval $(0,\infty)$ such that \cref{G1} holds
	and neither $h_1$ nor $h_2$ vanishes a.e.\ on some neighbourhood of the left endpoint $0$.
	Then the following statements are equivalent.
	\begin{Enumerate}
	\item
		${\displaystyle
			\liminf_{y\to\infty}\frac{\Im q_H(iy)}{|q_H(iy)|}=0
		}$.
	\item
		${\displaystyle
			\liminf_{t\to 0}\frac{\det M(t)}{m_1(t)m_2(t)}=0
		}$.
	\item
		For each $T\in(0,\infty)$ there exists a sequence $(s_n)_{n\in\bb N}$ with $s_n\to 0$, such that
		for all $\gamma\in[0,1)$, and all open intervals $I,J\subseteq\bb R\setminus\{0\}$ with
		$\ov I\cap\ov J=\emptyset$ and at least one of $I$ and $J$ being bounded, the following limit relations hold:
		\begin{align}
			& \lim_{n\to\infty}\Bigl[\lambda\Bigl(
			(0,T)\cap\mf t_{s_n}\bigl(\tfrac 1{s_n}\sigma_H^{-1}([0,\gamma])\bigr)
			\Bigr)\Bigr]=0,
			\label{G69}
			\\[1ex]
			&
			\lim_{n\to\infty}\Bigl[
			\lambda\Bigl((0,T)\cap\mf t_{s_n}\bigl(\pi_{H,s_n}^{-1}(I)\bigr)\Bigr)\cdot
			\lambda\Bigl((0,T)\cap\mf t_{s_n}\bigl(\pi_{H,s_n}^{-1}(J)\bigr)\Bigr)
			\Bigr]=0.
			\nonumber
		\end{align}
	\end{Enumerate}
\end{Theorem}

\noindent
Also in this case, the analogous addition holds.

\begin{GenericTheorem}{Addition to \Cref{G15}}
	Assume that, in addition to the assumptions of \Cref{G15}, relations \textup{\cref{G99}}
	and \textup{\cref{G100}} hold.
	Then the equivalent properties \Enumref{1}, \Enumref{2}, \Enumref{3} in \Cref{G15}
	are further equivalent to the following condition.
	\begin{Enumerate}
	\setcounter{enumi}{3}
	\item
		There exists a sequence $(t_n)_{n\in\bb N}$ with $t_n\to 0$, such that,
		for all $\gamma$ and $I,J$ as in \Cref{G15}\,\Enumref{3}, we have
		\begin{align}
			& \lim_{n\to\infty}\biggl[\frac 1{t_n}
			\lambda\Bigl((0,t_n)\cap\sigma_H^{-1}\bigl([0,\gamma]\bigr)\Bigr)\Bigr]=0,
			\label{G68}
			\\[1ex]
			& \lim_{n\to\infty}\biggl[
			\frac 1{t_n}\lambda\Bigl((0,{t_n})\cap\pi_H^{-1}(I)\Bigr)\cdot
			\frac 1{t_n}\lambda\Bigl((0,{t_n})\cap\pi_H^{-1}(J)\Bigr)
			\biggr]=0.
			\label{G101}
		\end{align}
	\end{Enumerate}
\end{GenericTheorem}

\noindent
The conditions \cref{G69} and \cref{G68} can be rewritten in integral form in the
very same way as before.  Namely, \cref{G69} as
\[
	\lim_{n\to\infty}\int_0^T
	\frac{\det H\bigl(s_n\mf t_{s_n}^{-1}(t)\bigr)}{(h_1h_2)\bigl(s_n\mf t_{s_n}^{-1}(t)\bigr)}
	\DD t=0,
\]
and \cref{G68} as
\[
	\lim_{n\to\infty}\frac 1{t_n}\int_0^{t_n}\frac{\det H(t)}{(h_1h_2)(t)}\DD t=0.
\]

\subsection*{Two examples}

Let us illustrate \Cref{G14,G15} with two examples. The first one demonstrates a standard situation; it will be revisited in a
more general form in \Cref{G98} of the present paper, and in the forthcoming paper \cite{langer.pruckner.woracek:asysupp}.
The second example demonstrates a more peculiar situation, where $\frac{\Im q_H}{|q_H|}$ oscillates.

\begin{Example}\label{G77}
	Let $\alpha_1,\alpha_2>0$, $\beta_1,\beta_2\in\bb R$, set
	\[
		\alpha_3\DE\frac{\alpha_1+\alpha_2}2,\quad \beta_3\DE\frac{\beta_1+\beta_2}2
		,
	\]
	and consider the Hamiltonian
	\[
		H(t)\DE
		\begin{pmatrix}
			t^{\alpha_1-1}|\log t|^{\beta_1} & t^{\alpha_3-1}|\log t|^{\beta_3}
			\\[2mm]
			t^{\alpha_3-1}|\log t|^{\beta_3} & t^{\alpha_2-1}|\log t|^{\beta_2}
		\end{pmatrix},
		\qquad t\in(0,\infty)
		.
	\]
	For this example a computation shows the following facts (this is elementary and we skip details):
	\begin{Enumerate}
	\item for $y\to\infty$,
		\[
			A_H(y)\asymp y^{\frac{\alpha_2-\alpha_1}{\alpha_1+\alpha_2}}
			(\log y)^{\frac{\beta_1\alpha_2-\beta_2\alpha_1}{\alpha_1+\alpha_2}}
			.
		\]
	\item We have
		\begin{align*}
			L_H(y)\asymp \Im q_H(iy)\asymp |q_H(iy)|\asymp A_H(y) & \quad \text{if}\ \alpha_1\neq\alpha_2
			\\[2mm]
			L_H(y)\lesssim \Im q_H(iy)\ll |q_H(iy)|\asymp A_H(y) & \quad \text{if}\ \alpha_1=\alpha_2
		\end{align*}
	\item The situation that $\lim_{y\to\infty}\frac{\Im q_H(iy)}{|q_H(iy)|}=0$, equivalently that
		$\lim_{t\to 0}\frac{\det M(t)}{m_1(t)m_2(t)}=0$, appears only when $q_H(iy)$ grows very
		slowly. In fact, if $\alpha_1=\alpha_2$, then
		\[
			A_H(y)\asymp(\log y)^{\frac{\beta_1-\beta_2}{2}},\quad
			\frac{L_H(y)}{A_H(y)}=\frac 1{(\log y)^2}
			.
		\]
	\end{Enumerate}
\end{Example}

\begin{Example}\label{G97}
	Let $(t_n)_{n\in\bb N}$ be a strictly decreasing sequence of positive numbers
	such that $\frac{t_{n+1}}{t_n}\to0$ (and hence $t_n\to0$),
	set $t_0\DE\infty$ and consider the partition $(0,\infty)=I_+\cup I_-$ where
	\[
		I_+ \DE \bigcup_{k=1}^\infty\bigl[t_{2k},t_{2k-1}\bigr),
		\qquad
		I_- \DE \bigcup_{k=0}^\infty\bigl[t_{2k+1},t_{2k}\bigr).
	\]
	Further, let $\varphi_+,\varphi_-\in(0,\pi)\setminus\{\frac\pi2\}$
	with $\varphi_+\ne\varphi_-$ and define the Hamiltonian $H$ by
	\[
		H(t) =
		\begin{pmatrix}
			\cos^2\varphi(t) & \cos\varphi(t)\cdot\sin\varphi(t)
			\\[1ex]
			\cos\varphi(t)\cdot\sin\varphi(t) & \sin^2\varphi(t)
		\end{pmatrix}
	\]
	where
	\[
		\varphi(t) =
		\begin{cases}
			\varphi_+, & t\in I_+,
			\\[1ex]
			\varphi_-, & t\in I_-.
		\end{cases}
	\]
	Clearly, \cref{G99} and \cref{G100} are satisfied, so that we can apply
	the Additions to \Cref{G14,G15}.
	Since $\sigma_H(t)=1$ for $t>0$, the limit relation \cref{G62},
	and hence also \cref{G68}, holds for every $\gamma\in[0,1)$.
	Let us now check whether \cref{G66} and \cref{G101} are satisfied.
	Since $\varphi(t)=\varphi_H(t)$, where $\varphi_H(t)$ is as in \cref{G58},
	it follows from \cref{G94} that
	\[
		\pi_H(t)=\sgn\bigl(\tfrac\pi2-\varphi(t)\bigr)\cdot\tan^2\varphi(t)
		= \sgn\bigl(\tfrac\pi2-\varphi_\pm\bigr)\cdot\tan^2\varphi_\pm
		\ED c_\pm
		\qquad\text{when} \ t\in I_\pm.
	\]
	The limit relations \cref{G66} and \cref{G101} hold trivially
	whenever $I\cap\{c_+,c_-\}=\emptyset$ or $J\cap\{c_+,c_-\}=\emptyset$.
	By symmetry, we only have to consider the case when $c_+\in I$ and $c_-\in J$,
	which we assume in the following.
	For $t>0$ we have
	\begin{align*}
		\lambda\bigl((0,t)\cap\pi_H^{-1}(I)\bigr)
		&= \lambda\bigl((0,t)\cap I_+\bigr)
		=
		\begin{cases}
			t-t_{2n}+\sum\limits_{k=n+1}^\infty\bigl(t_{2k-1}-t_{2k}\bigr),
			& t\in\bigl[t_{2n},t_{2n-1}\bigr),
			\\[1ex]
			\sum\limits_{k=n+1}^\infty\bigl(t_{2k-1}-t_{2k}\bigr),
			& t\in\bigl[t_{2n+1},t_{2n}\bigr),
		\end{cases}
		\\[1ex]
		\lambda\bigl((0,t)\cap\pi_H^{-1}(J)\bigr)
		&= \lambda\bigl((0,t)\cap I_-\bigr)
		=
		\begin{cases}
			\sum\limits_{k=n}^\infty\bigl(t_{2k}-t_{2k+1}\bigr),
			& t\in\bigl[t_{2n},t_{2n-1}\bigr),
			\\[1ex]
			t-t_{2n+1}+\sum\limits_{k=n+1}^\infty\bigl(t_{2k}-t_{2k+1}\bigr),
			& t\in\bigl[t_{2n+1},t_{2n}\bigr).
		\end{cases}
	\end{align*}
	Set $F(t)\DE\frac1t\lambda((0,t)\cap I_+)\cdot\frac1t\lambda((0,t)\cap I_-)$.
	Then
	\[
		F(t_{2n}) \le \frac{1}{t_{2n}}\sum_{k=n+1}^\infty\bigl(t_{2k-1}-t_{2k}\bigr)
		\le \frac{t_{2n+1}}{t_{2n}} \to 0
	\]
	as $n\to\infty$ and, similarly, $F(t_{2n+1})\to0$.
	This shows that \cref{G101} is satisfied and hence also \Enumref{1}
	in \Cref{G15}.
	On the other hand, for $n\in\bb N$ such that $\frac{t_{2n}}{t_{2n-1}}\le\frac12$,
	we have
	\[
		F(2t_{2n}) \ge \frac{1}{2t_{2n}}\bigl(2t_{2n}-t_{2n}\bigr)
		\cdot\frac{1}{2t_{2n}}\bigl(t_{2n}-t_{2n+1}\bigr)
		= \frac14\Bigl(1-\frac{t_{2n+1}}{t_{2n}}\Bigr)
		\to \frac14
	\]
	as $n\to\infty$.  This implies that \cref{G66} is not fulfilled
	and hence neither is \Enumref{1} in \Cref{G14}.
	To summarise, \Cref{G14,G15} show that
	\[
		\liminf_{y\to\infty}\frac{\Im q_H(iy)}{|q_H(iy)|}=0
		\qquad\text{and}\qquad
		\limsup_{y\to\infty}\frac{\Im q_H(iy)}{|q_H(iy)|}>0.
	\]
\end{Example}

\section{Preliminaries}
\label{G56}
\subsection[\textcolor{Dandelion}{Convergence of Hamiltonians}]{Convergence of Hamiltonians}
\label{G72}

We use the following notation for Hamiltonians on a finite or infinite interval.

\begin{Definition}\label{G46}
	Let $T\in(0,\infty]$.
	\begin{Enumerate}
	\item
		$\HamInt T$ is the set of all measurable functions $H\DF(0,T)\to\bb R^{2\times 2}$
		(up to equality a.e.)
		such that $H(t)\ge 0$ and $\tr H(t)>0$ a.e.;
	\item
		$\NHamInt T$ is the set of all $H\in\HamInt T$ such that $\tr H(t)=1$ a.e.;
	\item
		$\scHamInt T$ is the set of all $H\in\NHamInt T$ that are constant and satisfy
		$\det H(t)=0$ a.e.
	\end{Enumerate}
	If $T=\infty$, we often drop $T$ from the notation and just write
	$\Ham$, $\NHam$ and $\scHam$ instead of $\HamInt\infty$, $\NHamInt\infty$
	and $\scHamInt\infty$ respectively.
\end{Definition}

\noindent
We recall how $\NHam$ can be topologised appropriately.
This is already used in the work of L.~de~Branges.  An explicit formulation is given
in \cite{remling:2018};
for a more structural approach see \cite{pruckner.woracek:limp},
which we use as our main reference in the following.

For each $T<\infty$ the set $\NHamInt T$ is a subset of $L^1((0,T),\bb R^{2\times 2})$,
and hence naturally topologised with the $\|\Dummy\|_1$-topology or the weak $L^1$-topology.
It turns out that the latter is more suitable because the
weak $L^1$-topology on $\NHamInt T$ is compact and metrisable;
see \cite[Lemma~2.3]{pruckner.woracek:limp}.

Now consider the family $(\NHamInt T)_{T\in(0,\infty)}$ with the restriction maps
$\rho_T^{T'}\DF\NHamInt{T'}\to\NHamInt T$ for $T\le T'$.
The set $\NHam$ can be naturally viewed as the inverse limit of this family:
every function on $(0,\infty)$ can be identified with the family of all its restrictions
to finite intervals.
Endowed with the inverse limit topology (see, e.g.\ \cite[\S I.4.4]{bourbaki:1966}),
where we use the weak $L^1$-topology on $\NHamInt T$, the set $\NHam$ becomes a
compact metrisable space; see \cite[Lemma~2.9]{pruckner.woracek:limp}.
The map that assigns to a Hamiltonian $H$ its Weyl coefficient $q_H$ is continuous
when the set of Nevanlinna functions is endowed with the topology of locally uniform
convergence; see \cite[Theorem~2.12]{pruckner.woracek:limp}.

Throughout the remainder of the paper we often deal with limit points of
families of Hamiltonians.  In general, for a net $(x_i)_{i\in I}$ in some
topological space $X$, we denote by $\LP(x_i)_{i\in I}$ the set of its limit points, i.e.\
\[
	\LP(x_i)_{i\in I}\DE
	\bigl\{x\in X\DS
	\exists \text{ subnet }(x_{i(j)})_{j\in J}\DP\lim_{j\in J}x_{i(j)}=x\bigr\}.
\]
If there is a need to specify the topology, we shall add an index.
For example, if $X$ is a normed space, we write $\LP_{\|\Dummy\|}(x_i)_{i\in I}$
for limit points w.r.t.\ the norm topology, and $\LP_w(x_i)_{i\in I}$ for limit points
w.r.t.\ the weak topology.

\begin{Remark}\label{G27}
	In our context the space $X$ is usually metrisable, and the index set $I$
	is $\bb N$, $(0,1]$ or $[1,\infty)$, each endowed with the natural order
	(or the reverse order in the case of $(0,1]$).
	In these situations one can restrict attention to subsequences rather than subnets:
	\[
		\LP(x_i)_{i\in I}
		= \bigl\{x\in X\DS\exists \text{ subsequence }(x_{i_n})_{n\in\bb N}\DP\lim_{n\to\infty}x_{i_n}=x\bigr\}.
	\]
	Note that in the cases when $I=(0,1]$ or $I=[1,\infty)$, then $i_n\to0$
	or $i_n\to\infty$ respectively.
\end{Remark}

\noindent
We need the following simple fact about constant singular limit points.
It is proved using the compactness of $\NHam$, continuity of the restriction maps
$\rho_T\DF\NHam\to\NHamInt{T}$, and the obvious fact that
\begin{equation}\label{G26}
	\scHam=\bigl\{H\in\NHam\DS \forall T>0\DP \rho_T(H)\in\scHamInt T\bigr\}.
\end{equation}

\begin{lemma}\label{G25}
	Let $(H_i)_{i\in I}$ be a net in $\NHam$.
	Then the following two equivalences hold.
	\begin{Enumerate}
	\item
		${\displaystyle
			\LP(H_i)_{i\in I}\subseteq\scHam \quad\Leftrightarrow\quad
			\forall T>0\DP \LP_w(\rho_T(H_i))_{i\in I}\subseteq\scHamInt T
		}$
	\item
		${\displaystyle
			\LP(H_i)_{i\in I}\cap\scHam\neq\emptyset \quad\Leftrightarrow\quad
			\forall T>0\DP \LP_w(\rho_T(H_i))_{i\in I}\cap\scHamInt T\ne\emptyset
		}$
	\end{Enumerate}
\end{lemma}

\begin{proof}
	\hfill\\
	\Enumref{1}``$\Leftarrow$'':
	Assume that there exists $\mr H\in\LP(H_i)_{i\in I}\setminus\scHam$.
	By \cref{G26} we find $T>0$ such that $\rho_T(\mr H)\notin\scHamInt T$.
	Since $\rho_T$ is continuous, we have $\rho_T(\mr H)\in\LP_w(\rho_T(H_i))_{i\in I}$.
	\\
	\Enumref{1}``$\Rightarrow$'':
	Assume that there exist $T>0$ and
	$\wt H_T\in\LP_w(\rho_T(H_i))_{i\in I}\setminus\scHamInt T$.
	Since $\NHam$ is compact and $\rho_T$ is continuous, we find $\mr H\in\LP(H_i)_{i\in I}$
	such that $\rho_T(\mr H)=\wt H_T$. Clearly, $\mr H\notin\scHam$.
	\\
	\Enumref{2}``$\Rightarrow$'':
	Assume that there exists $\mr H\in\LP(H_i)_{i\in I}\cap\scHam$.
	Continuity of $\rho_T$ yields $\rho_T(\mr H)\in\LP(H_i)_{i\in I}\cap\scHamInt T$
	for all $T>0$.
	\\
	\Enumref{2}``$\Leftarrow$'':
	Assume that, for each $T>0$, there exists
	$\wt H_T\in\LP_w(\rho_T(H_i))_{i\in I}\cap\scHamInt T$.
	Since $\NHam$ is compact and $\rho_T$ is continuous, we find
	$\mr H_T\in\LP(H_i)_{i\in I}$ such that $\rho_T(\mr H_T)=\wt H_T$.
	Again by compactness, there exists a limit point $\mr H\in\LP(\mr H_T)_{T>0}$,
	say $\mr H=\lim_{n\to\infty}\mr H_{t_n}$ with some sequence $t_n\to\infty$.
	Then $\mr H\in\LP(H_i)_{i\in I}$, and for each $T>0$ we have
	\[
		\rho_T(\mr H)
		= \lim_{n\to\infty}{\mkern-8mu^w\,}\rho_T(\mr H_{t_n})
		= \lim_{\substack{n\to\infty \\[0.2ex] t_n\ge T}}{\mkern-8mu^w\,}
		\rho_T(\rho_{t_n}(\mr H_{t_n}))
		= \lim_{\substack{n\to\infty \\[0.2ex] t_n\ge T}}{\mkern-8mu^w\,}
		\underbrace{\rho_T(\wt H_{t_n})}_{\in\scHamInt T}\in\scHamInt T.
	\]
	For the last inclusion recall that $\scHamInt T$ is $\|\Dummy\|_1$-compact
	as a homeomorphic image of $\bb R\cup\{\infty\}$,
	see \cite[\S2.3]{pruckner.woracek:limp}, and hence also weakly closed.
	Again referring to \cref{G26} we obtain $\mr H\in\scHam$.
\end{proof}

\noindent
We also need the Weyl coefficients for constant Hamiltonians with zero determinant,
which can be found by an elementary calculation;
see \cite[Example~2.2(1)]{eckhardt.kostenko.teschl:2018}\footnote{
	In \cite{eckhardt.kostenko.teschl:2018}
	a different sign convention is used, namely the equation $y'(t)=-zJH(t)y(t)$
	is studied instead of \eqref{G60}.
	The corresponding Weyl coefficient is $\widetilde q_H(z)=-q_H(-z)$.
}.

\begin{lemma}\label{G28}
	Let $H$ as in \textup{\cref{G54}} be a constant Hamiltonian such that $h_3^2=h_1h_2$.
	Then
	\[
		q_H(z) =
		\begin{cases}
			\frac{h_3}{h_2} & \text{if} \ h_2\ne0,
			\\[1ex]
			\infty & \text{if} \ h_2=0.
		\end{cases}
	\]
\end{lemma}

\subsection[\textcolor{Dandelion}{Reparameterisation}]{Reparameterisation}
\label{G70}

Reparameterisation is the equivalence relation on the set of all Hamiltonians
defined as follows.

\begin{Definition}\label{G52}
	Two Hamiltonians $H$ and $\widehat H$, defined on respective intervals $[a,b)$
	and $[\hat a,\hat b)$, are called reparameterisations of each other
	if there exists a function $\varphi\DF[\hat a,\hat b)\to[a,b)$ that is
	strictly increasing, bijective and absolutely continuous with absolutely continuous
	inverse such that
	\begin{equation}\label{G53}
		\widehat H(x)=(H\circ\varphi)(x)\cdot\varphi'(x),\qquad x\in[\hat a,\hat b)\text{ a.e.}
	\end{equation}
\end{Definition}

\noindent
If $H$ and $\widehat H$ are related as in \cref{G53} and $y$ is a solution of \cref{G60},
then $y\circ\varphi$ is a solution of \cref{G53} with $H$ replaced by $\widehat H$.
Similarly, the fundamental solutions satisfy $\hatW(x,z)=W(\varphi(x),z)$
and hence
\begin{equation}\label{G115}
	q_{\widehat H}(z) = q_H(z).
\end{equation}
Moreover, the following obvious transformation rules hold:
\begin{alignat}{2}
	& \hatM=M\circ\varphi,\qquad
	&& \tr\widehat H(s)=\tr H(\varphi(s))\cdot \varphi'(s),
	\nonumber
	\\[1ex]
	& \pi_{\widehat H}=\pi_H\circ\varphi,\qquad
	&& \sigma_{\widehat H}=\sigma_H\circ\varphi.
	\label{G51}
\end{alignat}
Based on the transformation rule for the trace,
we see that every equivalence class of Hamiltonians modulo reparameterisation contains
exactly one element that is defined on the interval $(0,\infty)$ and is trace-normalised,
i.e.\ whose trace is equal to $1$ a.e.
In fact, given a Hamiltonian $H$ defined on some interval $(a,b)$,
we set $\mf t(t)\DE\int_a^t\tr H(x)\DD x$ and use $\varphi\DE\mf t^{-1}$.
This function is admissible to make a reparameterisation, since $\tr H(t)>0$ a.e.,
and hence $\mf t^{-1}$ is absolutely continuous.

Based on the transformation rule of the primitive $M$, we see that the quotient
in \cref{G64} transforms correspondingly.  Let us set
\begin{equation}\label{G37}
	d(H,t) \DE \frac{\det M(t)}{m_1(t)m_2(t)}.
\end{equation}
If $H$ and $\widehat H$ are related as in \cref{G53}, then
\begin{equation}\label{G82}
	d(\widehat H,s) = d\bigl(H,\varphi(s)\bigr).
\end{equation}

\subsection[\textcolor{Dandelion}{Hamiltonians starting with a vanishing diagonal entry}]{Hamiltonians starting with a vanishing
	diagonal entry}
\label{G71}

If a Hamiltonian starts with an interval where a diagonal entry vanishes,
then its Weyl coefficient has a simple, and extremal, asymptotics towards $+i\infty$.

Let $H$ be a Hamiltonian defined on some interval $(a,b)$.
Recall the following classical facts; see, e.g.\ \cite{kac.krein:1968a}.
\begin{Ilist}
\item
	Denote by $(a,\hat a)$ the maximal interval starting at $a$ such that $h_2(t)=0$
	for $t\in(a,\hat a)$ a.e., and assume that $\hat a>a$.  Then
	\[
		q_H(z) = \bigg(\int_a^{\hat a} h_1(t)\DD t\bigg)\cdot z + q_{H|_{(\hat a,b)}}(z).
	\]
	The leading order term is the term that is linear in $z$:
	\[
		\lim\limits_{y\to+\infty}\frac{1}{y}q_H(iy)=i \int_a^{\hat a}\tr H(t)\DD t,
	\]
	The case $\hat a = b$ is formally included and corresponds to $q_H \equiv \infty$.
\item
	Denote by $(a,\check a)$ the maximal interval starting at $a$ such that $h_1(t)=0$
	for $t\in(a,\check a)$ a.e., and assume that $\check a>a$.  Then
	\[
		q_H(z) = -\frac{1}{\bigl(\int_a^{\check a} h_2(t)\DD t\bigr)\cdot z - \frac{1}{q_{H|_{(\check a,b)}}(z)}\,}.
	\]
	Again the linear term gives the leading order asymptotics:
	\[
		\lim\limits_{y\to+\infty} y q_H(iy) = \frac{i}{\int_a^{\check a}\tr H(t)\DD t}.
	\]
	The case $\check a = b$ is formally included and corresponds to $q_H \equiv 0$.
\end{Ilist}
Translated to the spectral measure, $\hat a>a$ means that it should include a
``point mass at infinity'', and $\check a>a$ means that it has finite total mass.

In particular, the above relations show that, if $\hat a>a$ or $\check a>a$, then
\[
	\lim_{y\to\infty}\frac{\Im q_H(iy)}{|q_H(iy)|}=1.
\]

\subsection[\textcolor{Dandelion}{Representation of Hamiltonians by scalar functions}]{Representation of Hamiltonians by scalar functions}

We study the representation of a Hamiltonian $H$ by means of the functions $\sigma_H$
and $\zeta_H$, defined in \cref{G57} and \cref{G59} respectively,
a bit more systematically.
Denote by $\bb T$ the unit circle in the complex plane and, for $0<T\leq\infty$, set
\[
	\ms L(T)\DE
	\bigl\{f:(0,T)\to[0,1]\times\bb T\DS f\text{ measurable}\bigr\}\big/_\sim,
\]
where $f\sim g$ means that $f$ and $g$ coincide almost everywhere.
As usual, we suppress explicit notation of equivalence classes.
Moreover, we write a function $f\in\ms L(T)$ generically as a pair $f=(\sigma,\zeta)$
with $\sigma\DF(0,T)\to[0,1]$ and $\zeta\DF(0,T)\to\bb T$.

The set $\ms L(T)$ is contained in $L^1((0,T),\bb C^2)$ if $T$ is finite,
and in $L^1_{\Loc}([0,\infty),\bb C^2)$ if $T=\infty$.
In particular, for $T<\infty$, we may consider $\ms L(T)$ topologised with
the $\|\Dummy\|_1$-topology or the weak $L^1$-topology.

\begin{Definition}\label{G32}
	Let $0<T\le\infty$.  We define maps
	\begin{center}
	\begin{tikzcd}[column sep=normal]
		\ms L(T) \arrow[bend left=30]{r}{\Gamma} & \NHamInt T \arrow[bend left=30]{l}{\Xi}
	\end{tikzcd}
	\end{center}
	by
	\[
		\Gamma[\sigma,\zeta](t)\DE
		\frac 12
		\begin{pmatrix}
			1+\Re\zeta(t) & \sigma(t)\Im\zeta(t)
			\\[1ex]
			\sigma(t)\Im\zeta(t) & 1-\Re\zeta(t)
		\end{pmatrix}
	\]
	and $\Xi[H](t)\DE (\sigma_H(t),\zeta_H(t))$, where $\sigma_H$ and $\zeta_H$ are
	given by the formulae \cref{G57}, \cref{G58}, \cref{G59}.
\end{Definition}

\noindent
Let $(\sigma,\zeta)\in\ms L(T)$.
Introducing the rotation angle $\varphi\DF(0,T)\to[0,\pi)$ by $\zeta(t)=e^{2i\varphi(t)}$
we can rewrite
\[
	\Gamma[\sigma,\zeta]=
	\begin{pmatrix}
		1 & \sigma(t)
		\\
		\sigma(t) & 1
	\end{pmatrix}
	\odot
	\binom{\cos\varphi(t)}{\sin\varphi(t)}\binom{\cos\varphi(t)}{\sin\varphi(t)}^*.
\]
From this representation we see that $\Gamma$ is a left-inverse of $\Xi$:
given $H\in\NHamInt T$, the matrices $H$ and $\Gamma[\sigma_H,\zeta_H]$ both have trace $1$,
their quotients of diagonal entries coincide, and the relative size and sign of
their off-diagonal entries coincide.  Thus indeed
\[
	(\Gamma\circ\Xi)(H)=H.
\]
Furthermore, observe the following continuity property, which holds since $\ms L(T)$
is uniformly bounded: for each $T<\infty$ we find a constant $C>0$ such that,
for all $(\sigma_1,\zeta_1),(\sigma_2,\zeta_2)\in\ms L(T)$,
\begin{equation}\label{G33}
	\big\|\Gamma[\sigma_1,\zeta_1]-\Gamma[\sigma_2,\zeta_2]\big\|_1
	\le C\big\|(\sigma_1,\zeta_1)-(\sigma_2,\zeta_2)\big\|_1.
\end{equation}
In particular, for each $T<\infty$, the map $\Gamma\DF\ms L(T)\to\NHamInt T$
is $\|\Dummy\|_1$-to-$\|\Dummy\|_1$-continuous.

Note that the class of constant, singular, trace-normalised Hamiltonians can be
represented as follows:
\[
	\scHamInt T = \{\Gamma(1,\zeta):\zeta\in\bb T\}
\]
where we identify the constant $(1,\zeta)$ with the constant function in $\ms L(T)$.

\subsection{Nets with constant limit points}
\label{G74}

In the proof of the implication (iii)$\Rightarrow$(i) in \Cref{G14,G15}
we need the following fact about sequences in $L^1$-spaces which have only constant
limit points.  We do not know an explicit reference to the literature, and hence give
a complete proof.  In the formulation we tacitly identify $\bb C$ with
the $\mu$-a.e.\ constant functions in $L^1(\mu)$.

\begin{proposition}\label{G2}
	Let $\mu$ be a finite positive measure on a set $\Omega$ with $\mu\ne 0$, and
	let $(f_n)_{n\in\bb N}$ be a sequence in $L^\infty(\mu)$
	with $\sup_{n\in\bb N}\|f_n\|_\infty<\infty$.
	We consider $(f_n)_{n\in\bb N}$ as a sequence in $L^1(\mu)$.
	Then the following two statements are equivalent:
	\begin{align}
		& \forall (f_{n_k})_{k\in\bb N}\text{ subsequence of }(f_n)_{n\in\bb N}\DP
		\LP_{\|\Dummy\|_1}(f_{n_k})_{k\in\bb N}\cap\bb C\ne\emptyset;
		\label{G10}
		\\[1ex]
		& \forall A,B\subseteq\bb C\text{ compact, disjoint}\DP
		\lim_{n\to\infty}\Bigl[\mu\bigl(f_n^{-1}(A)\bigr)\cdot\mu\bigl(f_n^{-1}(B)\bigr)\Bigr]=0.
		\label{G11}
	\end{align}
	If the equivalent conditions \textup{\cref{G10}} and \textup{\cref{G11}} hold, then
	\begin{equation}
	\label{G43}
		\LP_w(f_n)_{n\in\bb N} = \LP_{\|\Dummy\|_1}(f_n)_{n\in\bb N}\subseteq\bb C.
	\end{equation}
\end{proposition}

\begin{proof}
	Let us first settle ``\cref{G10}$\Rightarrow$\cref{G11}$\wedge$\cref{G43}'',
	which is easy to see.
	Assume that \cref{G10} holds, and let $n_k\to\infty$.
	Then we find a further subsequence $(f_{n_{k(l)}})_{l\in\bb N}$
	such that $f_{n_{k(l)}}\stackrel{\|\Dummy\|_1}{\longrightarrow}g$ with some constant $g$.
	Since $\|\Dummy\|_1$-convergence implies convergence in measure, we have
	\begin{equation}\label{G36}
		\lim_{l\to\infty}\mu\bigl(\{x\in\Omega\DS |f_{n_{k(l)}}(x)-g|\ge\varepsilon\}\bigr)=0
	\end{equation}
	for every $\varepsilon>0$.
	Now consider two compact disjoint subsets $A,B$ of $\bb C$.
	Then the point $g$ can belong to at most one of $A$ and $B$.
	For definiteness, assume that $g\notin A$. Then the distance $\Dist(A,g)$ is positive,
	and
	\[
		f_n^{-1}(A)\subseteq\bigl\{x\in\Omega\DS |f_n(x)-g|\ge\Dist(A,g)\bigr\}.
	\]
	Relation \cref{G36} implies that
	\[
		\lim_{l\to\infty}\mu\bigl(f_{n_{k(l)}}^{-1}(A)\bigr)=0,
	\]
	and hence also the limit in \cref{G11} along the subsequence $(n_{k(l)})_{l\in\bb N}$
	is zero.
	Since we started with an arbitrary sequence $(n_k)_{k\in\bb N}$,
	the limit relation \cref{G11} follows.

	Now let $f\in\LP_w(f_n)_{n\in\bb N}$ and choose a subsequence $(f_{n_k})_{k\in\bb N}$
	with $f_{n_k}\stackrel w\to f$.
	Then we find a further subsequence $(f_{n_{k(l)}})_{l\in\bb N}$ and a constant $g$
	such that $f_{n_{k(l)}}\stackrel{\|\Dummy\|_1}{\longrightarrow}g$.
	It follows that $f=g\in\LP_{\|\Dummy\|_1}(f_n)_{n\in\bb N}$.  We have thus shown that
	\[
		\LP_w(f_n)_{n\in\bb N}\subseteq\LP_{\|\Dummy\|_1}(f_n)_{n\in\bb N}\cap\bb C,
	\]
	and this implies \cref{G43}.

	We come to the converse implication ``\cref{G11}$\Rightarrow$\cref{G10}''.
	Assume from now on that \cref{G11} holds.  Moreover, since \cref{G11} is inherited
	by subsequences, it is enough to prove \cref{G10} for the sequence $(f_n)_{n\in\bb N}$
	itself.  Further, let us set $M\DE\sup_{n\in\bb N}\|f_n\|_\infty$.

	There exist a subsequence $(n_k)_{k\in\bb N}$ and $a\in\bb R$ such that
	\begin{equation}\label{G131}
		\lim_{k\to\infty}\frac{1}{\mu(\Omega)}\int_\Omega \Re f_{n_k}(x)\DD\mu(x) = a.
	\end{equation}
	Let $\varepsilon>0$ be arbitrary and consider the compact, disjoint sets
	\[
		A = \bigl\{z\in\bb C\DF \Re z\ge a+\varepsilon \;\wedge\; |z|\le M\bigr\},
		\qquad
		B = \bigl\{z\in\bb C\DF \Re z\le a+\tfrac{\varepsilon}{2} \;\wedge\; |z|\le M\bigr\};
	\]
	by assumption, \cref{G11} holds with these sets.
	Suppose that there exist a subsequence $(k(l))_{l\in\bb N}$ such
	that $\lim_{l\to\infty}\mu\bigl(f_{n_{k(l)}}^{-1}(B)\bigr)=0$.  Then
	\begin{align*}
		\int_\Omega \Re f_{n_{k(l)}}(x)\DD\mu(x)
		&= \int_{f_{n_{k(l)}}^{-1}(B)} \Re f_{n_{k(l)}}(x)\DD\mu(x)
		+ \int_{\Omega\setminus f_{n_{k(l)}}^{-1}(B)} \Re f_{n_{k(l)}}(x)\DD\mu(x)
		\\[1ex]
		&\ge -M\mu\bigl(f_{n_{k(l)}}^{-1}(B)\bigr)
		+ \Bigl(a+\frac{\varepsilon}{2}\Bigr)\mu\bigl(\Omega\setminus f_{n_{k(l)}}^{-1}(B)\bigr)
		\\[1ex]
		&\to \Bigl(a+\frac{\varepsilon}{2}\Bigr)\mu(\Omega),
		\qquad l\to\infty,
	\end{align*}
	which is a contradiction to \cref{G131}.
	Therefore \cref{G11} implies that $\lim_{k\to\infty}\mu\bigl(f_{n_k}^{-1}(A)\bigr)=0$,
	which is equivalent to
	\begin{equation}\label{G132}
		\lim_{k\to\infty}\mu\bigl(\bigl\{x\in\Omega\DF
		\Re f_{n_k}(x)\ge a+\varepsilon\bigr\}\bigr) = 0.
	\end{equation}
	In a similar way one shows that
	\[
		\lim_{k\to\infty}\mu\bigl(\bigl\{x\in\Omega\DF
		\Re f_{n_k}(x)\le a-\varepsilon\bigr\}\bigr) = 0,
	\]
	which, together with \cref{G132}, implies that $\Re f_{n_k}\to a$ in measure.
	Since $\Re f_{n_k}$ is uniformly integrable, it follows that
	$\Re f_{n_k}\stackrel{\|\Dummy\|_1}{\longrightarrow}a$.
	In a completely analogous way one can find a subsequence such
	that $\Im f_{n_{k(l)}}\stackrel{\|\Dummy\|_1}{\longrightarrow}b$
	with some $b\in\bb R$.  This proves \cref{G10}.
\end{proof}

\medskip

\noindent
In the context of Hamiltonians on a finite interval, \Cref{G2} implies the following fact.

\begin{corollary}\label{G3}
	Let $T<\infty$ and $(H_n)_{n\in\bb N}$ be a sequence in $\NHamInt T$,
	and denote by $\lambda$ the Lebesgue measure on $(0,T)$.  Assume that
	\begin{Enumerate}
	\item
		${\displaystyle
		\forall \gamma\in[0,1)\DP
		\lim_{n\to\infty}\lambda\bigl(\sigma_{H_n}^{-1}([0,\gamma])\bigr)=0
		}$,
	\item
		${\displaystyle
		\forall A,B\subseteq\bb T\text{ closed, disjoint}\DP
		\lim_{n\to\infty}\Bigl[\lambda\bigl(\zeta_{H_n}^{-1}(A)\bigr)
		\cdot\lambda\bigl(\zeta_{H_n}^{-1}(B)\bigr)\Bigr]=0
		}$.
	\end{Enumerate}
	Then
	\[
		\LP_{\|\Dummy\|_1}(H_n)_{n\in\bb N} = \LP_w(H_n)_{n\in\bb N} \subseteq \scHamInt T.
	\]
\end{corollary}

\begin{proof}
	We have to show that
	\[
		\LP_w(H_n)_{n\in\bb N}\subseteq\LP_{\|\Dummy\|_1}(H_n)_{n\in\bb N}\cap\scHamInt T.
	\]
	The condition \Enumref{1} says that $\sigma_{H_n}\to 1$ in measure.
	Since $|\sigma_{H_n}(t)|\le1$ for a.e.\ $t$, $\sigma_{H_n}$ tends to $1$
	also w.r.t.\ $\|\Dummy\|_1$.
	Consider a subsequence $(H_{n_k})_{k\in\bb N}$ of $(H_n)_{n\in\bb N}$
	that converges weakly to some $\mr H\in\NHamInt T$.
	By \Enumref{2}, we can apply \Cref{G2} to the sequence $(\zeta_{H_{n_k}})_{k\in\bb N}$.
	This provides us with a constant $\zeta\in\bb T$ and a further
	subsequence $(\zeta_{H_{n_{k(l)}}})_{l\in\bb N}$
	such that $\zeta_{H_{n_{k(l)}}}\to\zeta$ w.r.t.\ $\|\Dummy\|_1$.
	Recalling \cref{G33} we see that
	\[
		\|H_{n_{k(l)}}-\Gamma(1,\zeta)\|_1
		= \big\|\Gamma(\sigma_{H_{n_{k(l)}}},\zeta_{H_{n_{k(l)}}})-\Gamma(1,\zeta)\big\|_1
		\to 0,
	\]
	and hence $\mr H=\Gamma(1,\zeta)$ and $\mr H\in\LP_{\|\Dummy\|_1}(H_n)_{n\in\bb N}$.
\end{proof}

\subsection[\textcolor{Dandelion}{Estimates for imaginary part and modulus of the Weyl coefficient}]{Estimates for imaginary part and modulus of the Weyl coefficient}
\label{G84}

In this subsection we recall lower and upper estimates for $\Im q_H$ and $|q_H|$
on the positive imaginary axis.  This result is a special instance of
\cite[Theorem~1.1]{langer.pruckner.woracek:heniest-arXiv} with $q=\frac14$
and $\vartheta=\frac\pi2$ there and is used, in particular, in the proof
of the implication (i)$\Rightarrow$(ii) in \Cref{G14};
the estimates for the modulus are also used in the proof of the implication
(iii)$\Rightarrow$(i).

\begin{Proposition}\label{G85}
	Let $H$ be a Hamiltonian defined on the interval $(0,\infty)$ such that \textup{\cref{G1}} holds
	and neither $h_1$ nor $h_2$ vanishes a.e.\ on some neighbourhood of the left endpoint $0$,
	and let $m_i$ be as in \textup{\cref{G54}} and $d(H,t)$ as in \textup{\cref{G37}}.
	For $r>0$, let $\mr t(r)\in(0,\infty)$ be the unique number that satisfies
	\begin{equation}\label{G88}
		(m_1m_2)\bigl(\mr t(r)\bigr)=\frac{1}{(8r)^2}\,.
	\end{equation}
	With
	\begin{equation}\label{G89}
		A_H(r) \DE \sqrt{\frac{m_1(\mr t(r))}{m_2(\mr t(r))}},
		\qquad
		L_H(r) \DE A_H(r)d\bigl(H,\mr t(r)\bigr)
	\end{equation}
	the inequalities
	\[
		\frac{1}{44}A_H(r) \le |q_H(ir)| \le 44A_H(r),
		\qquad
		\frac{1}{64}L_H(r) \le \Im q_H(ir) \le \frac{79}{2}A_H(r)
	\]
	hold for all $r>0$.
\end{Proposition}

\noindent
Note that the mapping $t\mapsto (m_1m_2)(t)$ is a strictly increasing
bijection from $(0,\infty)$ onto itself, and therefore $\mr t(r)$ is uniquely
defined via \cref{G88}. The mapping $r\mapsto\mr t(r)$ is a strictly decreasing
bijection from $(0,\infty)$ onto itself. It is the inverse of the function
\[
	\mr r(t) \DE \frac{1}{8\sqrt{(m_1m_2)(t)}\,}.
\]

\subsection[\textcolor{Dandelion}{A weighted rescaling transformation}]{A weighted rescaling transformation}
\label{G73}

In order to study the behaviour of $q_H$ towards $i\infty$, we use a weighted
rescaling transformation on the set of Hamiltonians.
This is a variant of Y.~Kasahara's rescaling trick invented in \cite{kasahara:1975}
for Krein strings, and also used in slightly different forms in
\cite{kasahara.watanabe:2010,eckhardt.kostenko.teschl:2018,pruckner.woracek:limp,langer.pruckner.woracek:asysupp}.
The main idea of the rescaling is to zoom into a neighbourhood of the left endpoint $0$
when $s$ in the following definition tends to $0$.

\begin{Definition}\label{G105}
	Let $g_1,g_2:(0,\infty)\to(0,\infty)$ be continuous
	such that $g_1(s),g_2(s)\to\infty$ as $s\to0$.
	Further, let $T\in(0,\infty]$ and set $g_3(s)\DE\sqrt{g_1(s)g_2(s)}$.
	For every $s>0$ define the map $\genResc s{}\DF\HamInt T\to\HamInt{\frac{1}{s}T}$
	by
	\[
		(\genResc sH)(t) \DE
		\begin{pmatrix}
			sg_1(s)h_1(st) & sg_3(s)h_3(st)
			\\[1ex]
			sg_3(s)h_3(st) & sg_2(s)h_2(st)
		\end{pmatrix},
		\qquad t\in\bigl(0,\tfrac1s T\bigr).
	\]
\end{Definition}

\noindent
In the following we shall use two special choices of $g_1,g_2$, namely
\begin{flalign}
	&\hspace*{5ex} \text{Situation 1:} \quad g_1(s) = \frac{1}{m_1(s)}, \quad g_2(s) = \frac{1}{m_2(s)}
	\label{G103}
	\\
	& \text{or} &
	\nonumber
	\\
	&\hspace*{5ex} \text{Situation 2:} \quad g_1(s) = g_2(s) = \frac{1}{s}
	\quad\text{and}\quad H \ \text{satisfies \cref{G99} and \cref{G100};}
	\label{G104}
\end{flalign}
in both cases $g_1,g_2$ satisfy the assumptions in \Cref{G105}.
The functions in \cref{G103} are used in the proof of \Cref{G14,G15};
the functions in \cref{G104} are used in the proof of the additions of these theorems.

In the following lemma we collect how the quantities defined in \cref{G57}--\cref{G76},
\cref{G87} and \cref{G37} are transformed.

\begin{lemma}\label{G106}
	Let $g_1,g_2$ be as in \Cref{G105} and $H\in\Ham$.
	Then
	\begin{align}
		& M\bigl(\genResc sH,t\bigr) = \int_0^t (\genResc sH)(x)\DD x =
		\begin{pmatrix}
			g_1(s)m_1(st) & g_3(s)m_3(st)
			\\[1ex]
			g_3(s)m_3(st) & g_2(s)m_2(st)
		\end{pmatrix},
		\label{G107}
		\\[1ex]
		& \sigma_{\genResc sH}(t) = \sigma_H(st),
		\qquad
		\pi_{\genResc sH}(t) = \frac{g_2(s)}{g_1(s)}\pi_H(st),
		\qquad
		d(\genResc sH,t) = d(H,st),
		\label{G108}
	\end{align}
	If, in addition, \textup{\cref{G1}} holds, then
	\begin{equation}\label{G109}
		q_{\genResc sH}(z) = \frac{g_3(s)}{g_2(s)}q_H\bigl(g_3(s)z\bigr).
	\end{equation}
\end{lemma}

\begin{proof}
	Relations \cref{G107} and the first two equalities in \cref{G108} follow easily
	from the definitions.
	Further, \cref{G107} implies the third equality in \cref{G108}.
	Finally, \cref{G109} follows
	from \cite[Lemma~2.7]{eckhardt.kostenko.teschl:2018}
\end{proof}

\noindent
If the functions $g_1,g_2$ are as in \cref{G103}, the relation \cref{G109} yields
\begin{equation}
\label{G122}
	q_{\genResc{\mr t(r)}H}\Big(\frac z8\Big) = \frac 1{A(r)}q_H(rz)
	.
\end{equation}
In the following lemma we prove an a priori estimate for the modulus of
the Weyl coefficient of $\genResc sH$ at a particular point, which is used in the proof of \Cref{G14,G15}.
This property follows from the choice of $g_1,g_2$ in \cref{G103} in the
general case or from the assumption \cref{G100} in the additions to
the main theorems.

\begin{lemma}\label{G35}
	Let $H\in\Ham$ such that \textup{\cref{G1}} holds, let $g_1,g_2,g_3$ be as in \Cref{G105},
	and assume that \textup{\cref{G103}} or \textup{\cref{G104}} is satisfied.  Then
	\[
		\Big|q_{\genResc sH}\Bigl(\frac i8\Bigr)\Big| \asymp 1, \qquad s\in(0,1].
	\]
\end{lemma}
\begin{proof}
	If $g_1,g_2$ are as in \cref{G103}, the assertion is clear from \cref{G122} and \Cref{G85}.

	Assume that \cref{G104} holds.
	Set $x_s\DE\mr t(\frac18)$ where $\mr t(r)$ is the unique number that
	satisfies \cref{G88} for $\genResc sH$ instead of $H$.  Then
	\[
		g_1(s)m_1(sx_s)g_2(s)m_2(sx_s) = 1.
	\]
	This is equivalent to $(m_1m_2)(sx_s)=s^2$. The latter relation implies
	that $sx_s\to0$ as $s\to0$. Assumptions \cref{G99} and \cref{G100} yield $m_i(t)\asymp t$, $i=1,2$
	and hence
	\[
		A_{\genResc sH}\Bigl(\frac18\Bigr)
		= \sqrt{\frac{m_1(sx_s)}{m_2(sx_s)}} \asymp 1, \qquad s\in(0,1].
	\]
	We obtain from \Cref{G85} that
	\[
		\Big|q_{\genResc sH}\Bigl(\frac i8\Bigr)\Big|
		\asymp A_{\genResc sH}\Bigl(\frac18\Bigr) \asymp 1,
		\qquad s\in(0,1].
	\]
\end{proof}

\medskip
\noindent
In the proof of \Cref{G14,G15} in \Cref{G83} we also need the trace of
the primitive of the rescaled Hamiltonian.
Let $g_1,g_2$ be as in \Cref{G105} and $H\in\Ham$.
For $s>0$ set
\begin{equation}\label{G111}
	\trResc s(t) \DE \int_0^t \tr(\genResc sH)(x)\DD x
	= g_1(s)m_1(st)+g_2(s)m_2(st),
	\qquad t\in(0,\infty).
\end{equation}
Since $\trResc s'(t) = sg_1(s)h_1(st)+sg_2(s)h_2(st)>0$ a.e.,
the function $\trResc s$ is strictly increasing.
If, in addition, $H$ is in the limit point case at $\infty$, then $\trResc s$
is a bijection from $(0,\infty)$ onto itself.
Note that for the choice \cref{G103} we have $\trResc s=\mf t_s$.

\section{Proof of ``(i)\,{\boldmath$\Leftrightarrow$}\,(ii)'' in \Cref{G14,G15}}
\label{G123}

We use the following fact which also plays a role later.

\begin{lemma}
\label{G110}
	Let $H_s$, $s>0$, be the trace-normalised reparameterisation of $\genResc sH$, i.e.\ the Hamiltonian that satisfies
	\begin{equation}\label{G96}
		(\genResc sH)(t) = H_s\bigl(\trResc s(t)\bigr)\cdot\trResc s'(t),
	\end{equation}
	where $\trResc s(t)$ is defined in \textup{\cref{G111}}.  Moreover, let $T\in(0,\infty)$. Then
	\begin{align}
		\lim_{t\to 0}d(H,t)=0 \quad&\Leftrightarrow\quad \lim_{s\to 0}d(H_s,T)=0,
		\label{G113}
		\\[1ex]
		\liminf_{t\to 0}d(H,t)=0 \quad&\Leftrightarrow\quad \liminf_{s\to 0}d(H_s,T)=0.
		\label{G114}
	\end{align}
\end{lemma}
\begin{proof}
	Let $T\in(0,\infty)$ be arbitrary.  By \cref{G82}, \cref{G108} and \cref{G111}
	we have
	\begin{equation}\label{G18}
		d(H_s,T) = d\bigl(\genResc sH,\trResc s^{-1}(T)\bigr)
		= d\bigl(H,s\trResc s^{-1}(T)\bigr).
	\end{equation}
	Set
	\begin{equation}\label{G80}
		u(s) \DE s\trResc s^{-1}(T).
	\end{equation}
	The explicit form of $\trResc s(t)$ and the continuity of $g_1$ and $g_2$
	show that the function $(s,t)\mapsto\trResc s(t)$ is continuous from $(0,\infty)^2$
	to $(0,\infty)$.
	Moreover, for every $s>0$, the mapping $t\mapsto\trResc s(t)$ is a
	homeomorphism from $(0,\infty)$ onto itself.
	By the implicit function theorem as in, e.g.\ \cite{kumagai:1980}, the function
	$s\mapsto\trResc s^{-1}(T)$, and with it also $s\mapsto u(s)$, is continuous.
	Moreover, we have
	\begin{equation}\label{G112}
		T = \trResc s\Bigl(\frac{u(s)}s\Bigr)
		= g_1(s)m_1\bigl(u(s)\bigr)+g_2(s)m_2\bigl(u(s)\bigr)
	\end{equation}
	for all $s$.  Since $g_i(s)\to\infty$ as $s\to0$ by assumption (see \Cref{G105}),
	it follows that $\lim_{s\to 0}u(s)=0$.
	The assertions now follow from \cref{G18}.
\end{proof}

\begin{proof}[Proof of ``\textup{(i)}$\Leftrightarrow$\textup{(ii)}'' in \Cref{G14,G15}]
	Let $H$ be as in the formulation of the theorems.
\begin{Elist}
\item
	The implications ``(i)$\Rightarrow$(ii)'' in \Cref{G14,G15} are a direct consequence of
	\cite[Theorem~1.1]{langer.pruckner.woracek:heniest-arXiv} in the form of \Cref{G85}
	since this result implies
	\[
		\frac{\Im q_H(ir)}{|q_H(ir)|}
		\ge \frac{\frac{1}{64}L_H(r)}{44A_H(r)}
		= \frac{1}{2816}d\bigl(H,\mr t(r)\bigr)
	\]
	for every $r>0$. It remains to recall that $\mr t$, defined via \cref{G88},
	is a strictly decreasing bijection from $(0,\infty)$ onto itself.
\item
	In this step we show that
	\begin{align}
	\label{G124}
		\lim_{r\to\infty}\frac{\Im q_H(ir)}{|q_H(ir)|}=0
		& \quad\Leftrightarrow\quad
		\LP(\mc A_s)_{s\in(0,1]}\subseteq\scHam
		\\[1ex]
	\label{G125}
		\liminf_{r\to\infty}\frac{\Im q_H(ir)}{|q_H(ir)|}=0
		& \quad\Leftrightarrow\quad
		\LP(\mc A_s)_{s\in(0,1]}\cap\scHam\ne\emptyset.
	\end{align}
	Let $r_n\to\infty$. Then we have the equivalences
	\[
		\lim_{n\to\infty}\frac{\Im q_H(ir_n)}{|q_H(ir_n)|}=0
		\quad\Leftrightarrow\quad
		\lim_{n\to\infty}\Im q_{\genResc{\mr t(r_n)}H}\Bigl(\frac i8\Bigr)=0
		\quad\Leftrightarrow\quad
		\LP(\mc A_{\mr t(r_n)})_{n\in\bb N}\subseteq\scHam
	\]
	The first one holds because of \cref{G122} and \Cref{G35}, and the second by the
	maximum principle and compactness of $\Ham$.
	Remembering that $\mr t$ is a decreasing bijection we obtain \cref{G124} and \cref{G125}.

\item
	To prove the implication ``(ii)$\Rightarrow$(i)'' in \Cref{G14},
	assume that $\lim_{t\to 0}d(H,t)=0$.
	By \Cref{G110} we have $\lim_{s\to 0}d(H_s,T)=0$ for all $T>0$.
	Since $\tr H_s=1$ a.e., it follows that also $\lim_{s\to 0}\det M_s(T)=0$ for all $T>0$
	where $M_s$ is the primitive of $H_s$.

	Let $\hat H\in\LP(\genResc sH)_{s\in(0,1]}$, and denote its primitive by $\hat M$.
	Then $\det\hat M(T)=0$ for all $T$. This means that the whole interval $(0,\infty)$
	is indivisible for $\hat H$, i.e.\ $\hat H\in\scHam$.
	Now \cref{G124} yields the required assertion.
\item
	For ``(ii)$\Rightarrow$(i)'' in \Cref{G15} assume that $\liminf_{t\to 0}d(H,t)=0$. Then for each $T>0$ we have
	$\liminf_{s\to 0}d(H_s,T)=0$ and, arguing as above, obtain a limit point $H_T\in\LP(\genResc sH)_{s\in(0,1]}$
	for which the interval $(0,T)$ is indivisible. Let $\phi_T\in[0,\pi)$ be the type of this indivisible interval.
	Choose a sequence $(T_n)_{n\in\bb N}$ such that $(\phi_{T_n})_{n\in\bb N}$ converges, say, $\phi_{T_n}\to\phi$.
	Then $(H_{T_n})_{n\in\bb N}$ converges to the Hamiltonian for which $(0,\infty)$ is indivisible of type $\phi$.
	Since $\LP(\genResc sH)_{s\in(0,1]}$ is closed, we can refer to \cref{G125} to finish the proof.
\end{Elist}
\end{proof}

\section{Bounds for the off-diagonal entries and the rotation}
\label{G55}

In this section we show that the relative size, $\sigma_H(t)$, of the off-diagonal entries
of a Hamiltonian and its rotation, $\zeta_H(t)$, can be estimated from above by $d(H,t)$;
recall that the latter is defined in \cref{G37}.
These estimates are used in the proof of the implication (ii)$\Rightarrow$(iii)
in \Cref{G14,G15}.

We start with an estimate for the off-diagonal entry.
As usual, $\lambda$ denotes the Lebesgue measure.

\begin{proposition}\label{G19}
	Let $H\in\NHam$, and assume that
	neither $h_1$ nor $h_2$ vanishes a.e.\ on some neighbourhood of the left endpoint $0$.
	For each $\gamma\in(0,1)$ we have
	\begin{equation}\label{G4}
		\forall t>0\DP
		\frac 1t\lambda\Bigl((0,t)\cap\sigma_H^{-1}\bigl([0,\gamma]\bigr)\Bigr)
		\le \frac{1}{1-\gamma^2}d(H,t).
	\end{equation}
\end{proposition}

\begin{proof}
	Throughout the proof we fix $\gamma\in(0,1)$ and $t>0$ and often suppress $t$ notationally.
	To shorten notation, we set $I_\gamma\DE(0,t)\cap\sigma_H^{-1}([0,\gamma])$
	and $I_\gamma'\DE(0,t)\setminus I_\gamma$.
	Further, set
	\[
		\xi_1 \DE \biggl[\int_{I_\gamma'}\!h_1(s)\DD s\biggr]^{\frac12}\!, \quad
		\xi_2 \DE \biggl[\int_{I_\gamma}\!h_1(s)\DD s\biggr]^{\frac12}\!, \quad
		\eta_1 \DE \biggl[\int_{I_\gamma'}\!h_2(s)\DD s\biggr]^{\frac12}\!, \quad
		\eta_2 \DE \biggl[\int_{I_\gamma}\!h_2(s)\DD s\biggr]^{\frac12}\!,
	\]
	define the vectors
	\[
		\upxi \DE \begin{pmatrix} \xi_1 \\[0.5ex] \xi_2 \end{pmatrix},
		\qquad
		\upeta \DE \begin{pmatrix} \eta_1 \\[0.5ex] \eta_2 \end{pmatrix},
	\]
	and let $\|\Dummy\|$ and $\langle\Dummy,\Dummy\rangle$ be the Euclidian norm and the
	inner product in $\bb R^2$.
	Since neither $h_1$ nor $h_2$ vanishes a.e.\ on $(0,t)$ by assumption,
	we have $\upxi\ne0$, $\upeta\ne0$, and we can write
	\[
		\upxi = \|\upxi\|\begin{pmatrix} \cos\theta_1 \\[0.5ex] \sin\theta_1 \end{pmatrix},
		\qquad
		\upeta = \|\upeta\|\begin{pmatrix} \cos\theta_2 \\[0.5ex] \sin\theta_2 \end{pmatrix}
	\]
	with $\theta_1,\theta_2\in\bigl[0,\frac\pi2\bigr]$.
	Moreover, we set $\theta\DE\max\{\theta_1,\theta_2\}$.

	The relation $I_\gamma=\bigl\{s\in(0,t)\DS |h_3(s)|\le\gamma\sqrt{h_1(s)h_2(s)}\bigr\}$
	implies that
	\begin{align*}
		|m_3(t)| &\le \int_{I_\gamma'}|h_3(s)|\DD s + \int_{I_\gamma}|h_3(s)|\DD s
		\\[1ex]
		&\le \int_{I_\gamma'}\sqrt{h_1(s)h_2(s)}\DD s
		+ \int_{I_\gamma}\gamma\sqrt{h_1(s)h_2(s)}\DD s
		\\[1ex]
		&\le \xi_1\eta_1+\gamma\xi_2\eta_2
		= \langle A\upxi,\upeta\rangle
	\end{align*}
	with $A=\bigl(\begin{smallmatrix} 1 & 0 \\ 0 & \gamma \end{smallmatrix}\bigr)$.
	For the diagonal terms in $M(t)$ we have
	\[
		m_1(t) = \int_{I_\gamma'}h_1(s)\DD s + \int_{I_\gamma}h_1(s)\DD s
		= \xi_1^2+\xi_2^2 = \|\upxi\|^2,
	\]
	and, similarly, $m_2(t)=\|\upeta\|^2$, which leads to
	\[
		\frac{m_3(t)^2}{m_1(t)m_2(t)}
		\le \frac{[\langle A\upxi,\upeta\rangle]^2}{\|\upxi\|^2\,\|\upeta\|^2}
		\le \min\biggl\{\frac{\|A\upxi\|^2}{\|\upxi\|^2},
		\frac{\|A\upeta\|^2}{\|\upeta\|^2}\biggr\}.
	\]
	The latter quotients can be rewritten as follows:
	\[
		\frac{\|A\upxi\|^2}{\|\upxi\|^2} = \frac{\xi_1^2+\gamma^2\xi_2^2}{\|\upxi\|^2}
		= \cos^2\theta_1+\gamma^2\sin^2\theta_1
		= 1-(1-\gamma^2)\sin^2\theta_1
	\]
	and analogously, $\|A\upeta\|^2/\|\upeta\|^2=1-(1-\gamma^2)\sin^2\theta_2$,
	which yields
	\begin{equation}\label{G127}
		\frac{m_3(t)^2}{m_1(t)m_2(t)}
		\le \min\bigl\{1-(1-\gamma^2)\sin^2\theta_1,1-(1-\gamma^2)\sin^2\theta_2\bigr\}
		= 1-(1-\gamma^2)\sin^2\theta.
	\end{equation}

	Since $H$ is trace-normalised, we have
	\[
		\xi_2^2+\eta_2^2 = \int_{I_\gamma}h_1(s)\DD s + \int_{I_\gamma}h_2(s)\DD s
		= \lambda(I_\gamma),
		\qquad
		\|\upxi\|^2+\|\upeta\|^2 = t,
	\]
	which implies that
	\begin{equation}\label{G128}
		\frac{\lambda(I_\gamma)}{t}
		= \frac{\|\upxi\|^2\sin^2\theta_1+\|\upeta\|^2\sin^2\theta_2}{\|\upxi\|^2+\|\upeta\|^2}
		\le \sin^2\theta.
	\end{equation}
	Combining \cref{G127} and \cref{G128} we obtain
	\[
		d(H,t) = 1-\frac{m_3(t)^2}{m_1(t)m_2(t)}
		\ge (1-\gamma^2)\sin^2\theta
		\ge (1-\gamma^2)\frac{\lambda(I_\gamma)}{t},
	\]
	which proves \cref{G4}.
\end{proof}

\noindent
Now we come to an estimate for the rotation of $H$.

\begin{Proposition}\label{G23}
	Let $H\in\NHam$, and assume that neither $h_1$ nor $h_2$ vanishes a.e.\
	on some neighbourhood of the left endpoint $0$.
	For each pair of closed, disjoint subsets $A,B\subseteq\bb T$ there exists
	a constant $c(A,B)>0$, which is independent of $H$, such that
	\[
		\forall t>0\DP
		\frac 1t\lambda\Bigl((0,t)\cap\zeta_H^{-1}(A)\Bigr)
		\cdot\frac 1t\lambda\Bigl((0,t)\cap\zeta_H^{-1}(B)\Bigr)
		\le c(A,B)\cdot d(H,t).
	\]
\end{Proposition}

\noindent
Heading towards the proof of this proposition, we present two lemmata.
The first one is an easy observation, which shows how information about the
Hamiltonian $H$ on an interval $I\subseteq(0,\infty)$ can be used to estimate $d(H,t)$.
In these two lemmata we use the following notation, which extends the notation
of the primitive to functions that may vanish on sets of positive measure:
for a Hamiltonian $H$, $I\subseteq(0,\infty)$ and $t>0$, set
\[
	M(H\mathds{1}_I,t) \equiv
	\begin{pmatrix} m_1(H\mathds{1}_I,t) & m_3(H\mathds{1}_I,t) \\[1ex]
	m_3(H\mathds{1}_I,t) & m_2(H\mathds{1}_I,t) \end{pmatrix}
	\DE \int_{I\cap(0,t)}H(s)\DD s.
\]

\begin{lemma}\label{G21}
	Let $H\in\NHam$ and $I\subseteq(0,\infty)$.
	For $t>0$, we have
	\[
		d(H,t) \ge \frac{\det M(H\mathds{1}_I,t)}{t^2}
	\]
\end{lemma}

\begin{proof}
	The fact that $H$ is positive semi-definite gives $M(H,t)\ge M(H\mathds{1}_I,t)\ge 0$,
	and in turn
	\[
		\det M(H,t) \ge \det M(H\mathds{1}_I,t) \ge 0.
	\]
	Together with $m_i(H,t)\le t$, which is a consequence of $\tr H=1$ a.e., we obtain
	\[
		d(H,t) = \frac{\det M(H,t)}{m_1(H,t)m_2(H,t)}
		\ge \frac{\det M(H\mathds{1}_I,t)}{m_1(H,t)m_2(H,t)}
		\ge \frac 1{t^2}\det M(H\mathds{1}_I,t).
	\]
\end{proof}

\noindent
The second lemma contains the crucial estimates.
For $\alpha,\beta\in\bb R$ with $\alpha\le\beta$ we denote the corresponding
arc on $\bb T$ by
\[
	A[\alpha,\beta] \DE \bigl\{\exp(it)\DS \alpha\le t\le\beta\bigr\}.
\]

\begin{lemma}\label{G22}
	The following estimates hold.
	\begin{Enumerate}	
	\item
		Let $\phi_0,\psi_0$ satisfy $0\leq\phi_0<\psi_0\leq\pi$ and set
		\\[-0.5ex]
		\begin{tabular}{@{\kern30pt}lr}
			\parbox{60mm}{$I_1\DE\zeta_H^{-1}\bigl({A[-\phi_0,\phi_0]}\bigr)$, \\[3mm]
			$I_2\DE\zeta_H^{-1}\bigl({A[\psi_0,2 \pi-\psi_0]}\bigr)$.}
			&
			\raisebox{-28pt}{\footnotesize
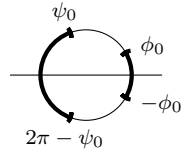
			\begin{tikzpicture}[scale=0.2]
				\draw (0,3) -- (10,3);
				\draw[thin] (5,3) circle (3);
				\draw[ultra thick] (8,3) arc(0:30:3)
					arc(120:140:1) arc(140:100:1) +(0,0.3) node[right]{$\phi_0$};
				\draw[ultra thick] (8,3) arc(0:-30:3)
					arc(-120:-140:1) arc(-140:-100:1) +(0,-0.3) node[right]{$-\phi_0$};
				\draw[ultra thick] (2,3) arc(180:250:3)
					arc(-20:0:1) arc(0:-40:1) +(-0.2,-1.1) node{$2\pi-\psi_0$};
				\draw[ultra thick] (2,3) arc(180:110:3)
					arc(20:0:1) arc(0:40:1) +(-0.2,1.1) node{$\psi_0$};
			\end{tikzpicture}
			}
		\end{tabular}
		\\[-0.5ex]
		Then, for all $H$ and $t>0$, we have
		\[
			d(H,t) \ge \sin^2\Bigl(\frac{\psi_0-\phi_0}{2}\Bigr)\cdot
			\frac 1t\lambda\bigl(I_1\cap(0,t)\bigr)
			\cdot\frac 1t\lambda\bigl(I_2\cap(0,t)\bigr).
		\]
	\item
		Let $\alpha,\beta\in(0,\pi]$ and set
		\\[1ex]
		\begin{tabular}{@{\kern30pt}lr}
			\parbox{55mm}{$I_1\DE\zeta_H^{-1}\bigl({A[0,\pi-\alpha]}\bigr)$, \\[3mm]
			$I_2\DE\zeta_H^{-1}\bigl({A[\pi,2\pi-\beta]}\bigr)$.}
			&
			\raisebox{-28pt}{\footnotesize
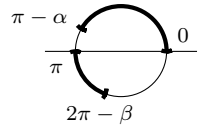
			\begin{tikzpicture}[scale=0.2]
				\draw (0,3) -- (10,3);
				\draw[thin] (5,3) circle (3);
				\draw[ultra thick] (8,3) arc(0:150:3)
					arc(240:260:1) arc(260:220:1) +(-0.0,0.5) node[left]{$\pi-\alpha$};
				\draw[ultra thick] (8,3)
					arc(270:290:1) arc(290:250:1) +(0.4,1) node[right]{$0$};
				\draw[ultra thick] (2,3) arc(180:250:3)
					arc(-20:0:1) arc(0:-40:1) +(-0.2,-1.1) node{$2\pi-\beta$};
				\draw[ultra thick] (2,3)
					arc(90:70:1) arc(70:110:1) +(0.2,-1.1) node[left]{$\pi$};
			\end{tikzpicture}
			}
		\end{tabular}
		\\[-1ex]
		Then, for all $H$ and $t>0$, we have
		\[
			d(H,t) \ge \sin^2\Bigl(\frac{\alpha}{2}\Bigr)\sin^2\Bigl(\frac{\beta}{2}\Bigr)
			\cdot\frac 1t\lambda\bigl(I_1\cap(0,t)\bigr)
			\cdot\frac 1t\lambda\bigl(I_2\cap(0,t)\bigr).
		\]
		The same holds for
		\\[-1ex]
		\begin{tabular}{@{\kern30pt}lr}
			\parbox{56mm}{$I_1\DE\zeta_H^{-1}\bigl({A[\pi+\alpha,2\pi]}\bigr)$, \\[3mm]
			$I_2\DE\zeta_H^{-1}\bigl({A[\beta,\pi]}\bigr)$.}
			&
			\raisebox{-17pt}{\footnotesize
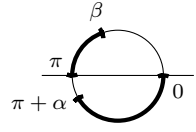
			\begin{tikzpicture}[scale=0.2]
				\draw (0,3) -- (10,3);
				\draw[thin] (5,3) circle (3);
				\draw[ultra thick] (8,3) arc(0:-150:3)
					arc(120:140:1) arc(140:100:1) +(-0.4,-0.5) node[left]{$\pi+\alpha$};
				\draw[ultra thick] (8,3)
					arc(90:110:1) arc(110:70:1) +(-0.4,-1) node[right]{$0$};
				\draw[ultra thick] (2,3) arc(180:110:3)
					arc(20:0:1) arc(0:40:1) +(-0.2,1.1) node{$\beta$};
				\draw[ultra thick] (2,3)
					arc(270:250:1) arc(250:290:1) +(-0.3,0.8) node[left]{$\pi$};
			\end{tikzpicture}
			}
		\end{tabular}
		\\[-1ex]
	\item
		Let $\alpha,\beta\in(0,\pi]$ satisfy $\alpha+\beta\leq\pi$ and set
		\\[-0.5ex]
		\begin{tabular}{@{\kern30pt}lr}
			\parbox{56mm}{$I_1\DE\zeta_H^{-1}\bigl({A[\beta,\pi-\alpha]}\bigr)$, \\[3mm]
			$I_2\DE\zeta_H^{-1}\bigl({A[\pi,2\pi]}\bigr)$.}
			&
			\raisebox{-17pt}{\footnotesize
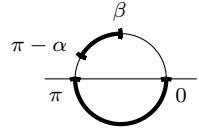
			\begin{tikzpicture}[scale=0.2]
				\draw (0,3) -- (10,3);
				\draw[thin] (5,3) circle (3);
				\draw[ultra thick] (8,3) arc(0:-180:3)
					arc(90:70:1) arc(70:110:1) +(0.2,-1.1) node[left]{$\pi$};
				\draw[ultra thick] (8,3)
					arc(90:110:1) arc(110:70:1) +(-0.4,-1) node[right]{$0$};
				\draw[ultra thick] (5,6) arc(90:150:3)
					arc(240:260:1) arc(260:220:1) +(-0.0,0.5) node[left]{$\pi-\alpha$};
				\draw[ultra thick] (5,6)
					arc(0:-20:1) arc(-20:20:1) +(0,1.1) node{$\beta$};
			\end{tikzpicture}
			}
		\end{tabular}
		\\[2mm]
		Then, for all $H$ and $t>0$, we have
		\[
			d(H,t) \ge \sin^2\Bigl(\frac{\min\{\alpha, \beta\}}{2}\Bigr)
			\cdot\frac 1t\lambda\bigl(I_1\cap(0,t)\bigr)
			\cdot\frac 1t\lambda\bigl(I_2\cap(0,t)\bigr).
		\]
		The same holds for
		\\[0.5ex]
		\begin{tabular}{@{\kern30pt}lr}
			\parbox{57mm}{$I_1\DE\zeta_H^{-1}\bigl({A[\pi+\alpha,2\pi-\beta]}\bigr)$, \\[3mm]
			$I_2\DE\zeta_H^{-1}\bigl({A[0,\pi]}\bigr)$.}
			&
			\raisebox{-27pt}{\footnotesize
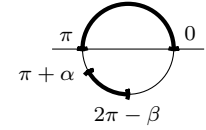
			\begin{tikzpicture}[scale=0.2]
				\draw (0,3) -- (10,3);
				\draw[thin] (5,3) circle (3);
				\draw[ultra thick] (8,3) arc(0:180:3)
					arc(270:250:1) arc(250:290:1) +(-0.3,0.8) node[left]{$\pi$};
				\draw[ultra thick] (8,3)
					arc(270:290:1) arc(290:250:1) +(0.4,1) node[right]{$0$};
				\draw[ultra thick] (5,0) arc(-90:-150:3)
					arc(-240:-260:1) arc(-260:-220:1) +(0.0,0) node[left]{$\pi+\alpha$};
				\draw[ultra thick] (5,0)
					arc(0:20:1) arc(20:-20:1) +(0,-1.1) node{$2\pi-\beta$};
			\end{tikzpicture}
			}
		\end{tabular}
		\\[-2.5ex]
	\end{Enumerate}
\end{lemma}

\begin{proof}
	\hfill
\begin{Elist}
\item We start with a general calculation.
	Let $K_1,K_2\subseteq [0,t]$ be disjoint and set $K\DE K_1 \cup K_2$.
	We can use the inequality $|h_3| \le \sqrt{h_1 h_2}$ and the Cauchy--Schwarz inequality
	in the last step to obtain
	\begin{align}
		\det M(H\mathds{1}_K,t)
		&= m_1( H \mathds{1}_K,t)m_2( H \mathds{1}_K,t) - m_3^2( H \mathds{1}_K,t)
		\nonumber
		\\[1ex]
		&= \biggl(\int_{K_1}h_1(x)\DD x +\int_{K_2}h_1(x)\DD x\biggr)
		\biggl(\int_{K_1}h_2(x)\DD x +\int_{K_2}h_2(x)\DD x\biggr)
		\nonumber
		\\[0.5ex]
		&\quad - \biggl(\int_{K_1}h_3(x)\DD x +\int_{K_2}h_3(x)\DD x\biggr)^2
		\nonumber\displaybreak[0]
		\\[1ex]
		&\ge \int_{K_1}h_1(x)\DD x \int_{K_1}h_2(x)\DD x
		-\biggl(\int_{K_1}\sqrt{h_1(x)h_2(x)}\DD x\biggr)^2
		\nonumber
		\\[0.5ex]
		&\quad + \int_{K_2}h_1(x)\DD x \int_{K_2}h_2(x)\DD x
		-\left( \int_{K_2}\sqrt{h_1(x)h_2(x)}\DD x\right)^2
		\nonumber
		\\[0.5ex]
		&\quad + \int_{K_1}h_1(x)\DD x \int_{K_2}h_2(x)\DD x + \int_{K_2}h_1(x)\DD x \int_{K_1}h_2(x)\DD x
		\nonumber
		\\[0.5ex]
		&\quad - 2 \int_{K_1}h_3(x)\DD x \int_{K_2}h_3(x)\DD x
		\nonumber\displaybreak[0]
		\\[1ex]
		&\ge \int_{K_1}h_1(x)\DD x \int_{K_2}h_2(x)\DD x
		+ \int_{K_2}h_1(x)\DD x \int_{K_1}h_2(x)\DD x
		\nonumber
		\\[0.5ex]
		&\quad - 2 \int_{K_1}h_3(x)\DD x \int_{K_2}h_3(x)\DD x.
		\label{G17}
	\end{align}
	Using once more $|h_3|\le \sqrt{h_1 h_2}$ and the Cauchy--Schwarz inequality
	we arrive at a complete square:
	\begin{align}
		\det M(H \mathds{1}_K,t)
		&\ge \int_{K_1}h_1(x)\DD x \int_{K_2}h_2(x)\DD x
		+ \int_{K_2}h_1(x)\DD x \int_{K_1}h_2(x)\DD x
		\nonumber
		\\[0.5ex]
		&\quad -2\biggl[\int_{K_1}h_1(x)\DD x \int_{K_1}h_2(x)\DD x
		\int_{K_2}h_1(x)\DD x \int_{K_2}h_2(x)\DD x\biggr]^{\frac 12}
		\nonumber
		\\[1ex]
		&= \Biggl[\biggl(\int_{K_1}h_1(x)\DD x\int_{K_2}h_2(x)\DD x\biggr)^{\frac12}
		-\biggl(\int_{K_2}h_1(x)\DD x\int_{K_1}h_2(x)\DD x\biggr)^{\frac12}\Biggr]^2.
		\label{G20}
	\end{align}
\item
	Let $\alpha,\beta\in\bb R$ with $\alpha\le\beta$ and
	set $J\DE\zeta_H^{-1}(A[\alpha,\beta])$.  Then
	\begin{align}
		x\in J \quad&\Leftrightarrow\quad \zeta_H(x)\in A[\alpha,\beta]
		\nonumber
		\\[1ex]
		&\Leftrightarrow\quad \exists\, n\in\bb Z\DP
		\varphi_H(x)-n\pi\in\bigl[\tfrac\alpha2,\tfrac\beta2\bigr]
		\nonumber
		\\[1ex]
		&\Rightarrow\quad \exists\,\varphi\in\bigl[\tfrac\alpha2,\tfrac\beta2\bigr]\DP
		h_1(x)=\cos^2\varphi, \quad h_2(x)=\sin^2\varphi.
		\label{G91}
	\end{align}
	For the rest of the proof set $K_i\DE I_i\cap(0,t)$ for $i=1,2$,
	and $K\DE K_1\cup K_2$.  We consider the three cases in \Enumref{1}, \Enumref{2},
	\Enumref{3} separately.
\item
	Let us first consider the situation in item \Enumref{1}.
	It follows from \cref{G91} that
	\begin{alignat*}{3}	
		h_1(x) &\ge \cos^2\Bigl(\frac{\phi_0}{2}\Bigr), \quad &
		h_2(x) &\le \sin^2\Bigl(\frac{\phi_0}{2}\Bigr), \qquad && x\in K_1,
		\\[1ex]
		h_1(x) &\le \cos^2\Bigl(\frac{\psi_0}{2}\Bigr), \quad &
		h_2(x) &\ge \sin^2\Bigl(\frac{\psi_0}{2}\Bigr), \qquad && x\in K_2.
	\end{alignat*}
	This, together with \Cref{G21} and \cref{G20}, implies
	\begin{align*}
		d(H,t) &\ge \frac{1}{t^2}\det M(H\mathds{1}_K,t)
		\\[1ex]
		&\ge \frac{1}{t^2}\Biggl[\sqrt{\cos^2\Bigl(\frac{\phi_0}{2}\Bigr)\lambda(K_1)
		\sin^2\Bigl(\frac{\psi_0}{2}\Bigr)\lambda(K_2)}
		- \sqrt{\cos^2\Bigl(\frac{\psi_0}{2}\Bigr)\lambda(K_2)
		\sin^2\Bigl(\frac{\phi_0}{2}\Bigr)\lambda(K_1)}\,\Biggr]^2
		\\[1ex]
		&= \biggl[\cos\Bigl(\frac{\phi_0}{2}\Bigr)\sin\Bigl(\frac{\psi_0}{2}\Bigr)
		- \cos\Bigl(\frac{\psi_0}{2}\Bigr)\sin\Bigl(\frac{\phi_0}{2}\Bigr)\biggr]^2
		\cdot\frac{\lambda(K_1)\lambda(K_2)}{t^2}
		\\[1ex]
		&= \sin^2\Bigl(\frac{\psi_0-\phi_0}{2}\Bigr)\cdot\frac{\lambda(K_1)}{t}
		\cdot\frac{\lambda(K_2)}{t},
	\end{align*}
	which is the asserted statement in \Enumref{1}.
\item
	Next, we consider the situation in item \Enumref{2}.
	Here $h_3$ is non-negative on $I_1$ and non-positive on $I_2$, or vice versa.
	Thus, \Cref{G21} and inequality \cref{G17} yield
	\begin{align}
		d(H,t) &\ge \frac{1}{t^2}\det M(H \mathds{1}_K,t)
		\nonumber
		\\[1ex]
		&\ge \frac{1}{t^2}\biggl[\int_{K_1}h_1(x)\DD x \int_{K_2}h_2(x)\DD x
		+ \int_{K_2}h_1(x)\DD x \int_{K_1}h_2(x)\DD x\biggr].
		\label{G38}
	\end{align}
	By \cref{G91} we have the estimates $h_1(x)\ge\cos^2((\pi\pm\alpha)/2)=\sin^2(\alpha/2)$
	for $x\in K_1$, and $h_2(x) \ge \sin^2(\beta/2)$ for $x\in K_2$,
	and hence
	\[
		d(H,t) \ge \sin^2\Bigl(\frac{\alpha}{2}\Bigr)\sin^2\Bigl(\frac{\beta}{2}\Bigr)
		\cdot\frac{\lambda(K_1)}{t}\cdot\frac{\lambda(K_2)}{t}.
	\]
\item
	Finally, we consider the situation in item \Enumref{3}.
	Again $h_3$ is non-negative on $I_1$ and non-positive on $I_2$, or vice versa,
	and therefore we obtain \cref{G38}.
	Further, for $x\in I_1$ we have the estimates
	$h_1(x)\ge\cos^2((\pi\pm\alpha)/2)=\sin^2(\alpha/2)$ and $h_2(x)\ge\sin^2(\beta/2)$
	by \cref{G91}.  With \Cref{G21} we obtain
	\begin{align*}
		d(H,t) &\ge \sin^2\Bigl(\frac{\alpha}{2}\Bigr)\cdot\frac{\lambda(K_1)}{t^2}
		\int_{K_2}h_2(x)\DD x
		+ \sin^2\Bigl(\frac{\beta}{2}\Bigr)\cdot\frac{\lambda(K_1)}{t^2}
		\int_{K_2}h_1(x)\DD x
		\\[1ex]
		&\ge \min\biggl\{\sin^2\Bigl(\frac{\alpha}{2}\Bigr),
		\sin^2\Bigl(\frac{\beta}{2}\Bigr)\biggr\}
		\cdot\frac{\lambda(K_1)}{t^2}\int_{K_2}\bigl(h_1(x)+h_2(x)\bigr)\DD x
		\\[1ex]
		&= \sin^2\Bigl(\frac{\min\{\alpha,\beta\}}{2}\Bigr)
		\cdot\frac{\lambda(K_1)}{t}\cdot\frac{\lambda(K_2)}{t}.
	\end{align*}
\end{Elist}
\vspace*{-2ex}
\end{proof}

\begin{proof}[Proof of \Cref{G23}]
	Let $\dd$ denote the intrinsic metric on $\bb T$ which assigns to a pair of points
	the length of the shortest arc connecting them.
\begin{Elist}
\item
	As a first step we settle the case when $A,B\subseteq \bb T$ are two
	closed, disjoint arcs with lengths strictly less than $\pi$.
	Set $\delta\DE\dd(A,B)$, fix $t>0$, and set
	\[
		\nu_A \DE \frac{1}{t}\lambda\bigl(\zeta_H^{-1}(A)\cap(0,t)\bigr),
		\qquad
		\nu_B \DE \frac{1}{t}\lambda\bigl(\zeta_H^{-1}(B)\cap(0,t)\bigr).
	\]
	Based on \Cref{G22} we are going to show that
	\begin{equation}\label{G39}
		d(H,t) \ge \frac{\sin^4(\delta/4)}4\cdot\nu_A\nu_B.
	\end{equation}
	To this end, we distinguish four cases.
	\begin{Ilist}
	\item
		Assume that one of $A$ and $B$ is contained in $A[0,\pi]$, the other one
		is contained in $A[\pi,2\pi]$,
		and either $\dd(A,1)\le\dd(B,1)$ and $\dd(B,-1)\le\dd(A,-1)$,
		or $\dd(B,1)\le\dd(A,1)$ and $\dd(A,-1)\le\dd(B,-1)$.
		Then \Cref{G22}\,\Enumref{2} with the choice $\alpha=\beta=\delta/2$ yields
		\begin{equation}\label{G40}
			d(H,t)\geq\sin^4(\delta/4)\cdot\nu_A\nu_B,
		\end{equation}
		which is even stronger than \cref{G39}.
	\item
		Assume that one of $A$ and $B$ is contained in $A[0,\pi]$, the other one
		is contained in $A[\pi,2\pi]$,
		and either $\dd(A,1)\le\dd(B,1)$ and $\dd(A,-1)\le\dd(B,-1)$,
		or $\dd(B,1)\le\dd(A,1)$ and $\dd(B,-1)\le\dd(A,-1)$.
		Then \Cref{G22}\,\Enumref{3} with the choice $\alpha=\beta=\delta/2$
		yields $d(H,t)\geq\sin^2(\delta/4)\nu_A \nu_B$, which implies \cref{G40}.
	\item
		Assume that both $A$ and $B$ are contained $A[0,\pi]$, or both are
		contained in $A[\pi,2\pi]$.
		Then \Cref{G22}\,\Enumref{1} yields $d(H,t)\ge\sin^2(\delta/2)\nu_A\nu_B$,
		which implies \cref{G40}.
	\item
		Assume that neither of the above three cases takes place, and set
		\begin{alignat*}{2}
			A_1 &\DE A\cap A[0,\pi],\qquad & A_2 &\DE A\cap A[\pi,2\pi],
			\\[1ex]
			B_1 &\DE B\cap A[0,\pi],\qquad & B_2 &\DE B\cap A[\pi,2\pi].
		\end{alignat*}
		Then $A_i$ and $B_j$ are contained in $A[0,\pi]$ or $A[\pi,2\pi]$, and
		satisfy $\dd(A_i,B_j)\geq\delta$.
		Moreover, since the lengths of $A$ and $B$ are strictly less than $\pi$,
		the sets $A_i$ and $B_j$ are again closed arcs.
		The already settled cases can be applied to $A_i$ and $B_j$, which yields
		\[
			d(H,t) \ge \sin^4(\delta/4)\cdot\nu_{A_i}\nu_{B_j},\qquad i,j\in\{1,2\};
		\]
		cf.\ \cref{G40}.
		There is at least one choice of $i,j\in\{1,2\}$ such that $\nu_{A_i}\ge\nu_A/2$
		and $\nu_{B_j}\geq\nu_B/2$.
		Using this choice we obtain \cref{G39}.
	\end{Ilist}
\item
	The general case, namely when $A$ and $B$ are arbitrary closed, disjoint subsets
	of $\bb T$, is deduced by appropriately covering $A$ and $B$ with arcs.

	Set $\delta\DE\dd(A,B)$ and consider the open cover of $\bb T$ consisting of
	all open arcs with length $\delta/3$.  Since $A$ is compact, there exist
	finitely many of these arcs whose union covers $A$, say $A_1,\ldots,A_N$.
	In addition, we may assume that $A\cap A_i\neq\emptyset$ for all $i\in\{1,\ldots,N\}$.
	In the same way we obtain arcs $B_1,\ldots,B_{N'}$ whose union covers $B$ and
	such that each of them intersects $B$.

	We have $\dd\bigl(\ov{A_i},\ov{B_k}\bigr)\ge\delta/3$ for all $i$ and $k$ by construction,
	and \cref{G39} tells us that
	\[
		d(H,t) \ge \frac{\sin^4(\delta/12)}4\cdot\nu_{\ov{A_i}}\nu_{\ov{B_k}}.
	\]
	For each $t>0$ there is at least one choice of $i\in\{1,\ldots,N\}$
	and $k\in\{1,\ldots,N'\}$ such that $\nu_{\ov{A_i}}\ge\nu_A/N$
	and $\nu_{\ov{B_k}}\ge\nu_B/N'$.  Using this choice we arrive at
	\[
		d(H,t) \ge \frac{\sin^4(\delta/12)}{4NN'}\cdot\nu_A\nu_B.
	\]
	Note that the constants $\delta,N,N'$ only depend on $A$ and $B$, but not on $H$ or $t$.
\end{Elist}
\vspace*{-2ex}
\end{proof}

\section{Proof of equivalence with (iii) in \Cref{G14,G15}}
\label{G83}

We have now collected all necessary tools to carry out the proof of equivalence with 
condition (iii) in our main theorems.
Our plan to proceed is to first work with a modified variant of \Enumref{3},
namely ``(iii$'$)'' stated below,
and prove that ``\Enumref{2}$\Rightarrow$(iii$'$)$\Rightarrow$\Enumref{1}''.
After that we show ``(iii$'$)$\Leftrightarrow$\Enumref{3}'', which is elementary.

In the following we consider the weighted rescalings $\genResc sH$ of $H$
from \Cref{G105}.
For most part of the proof, $g_1$ and $g_2$ are arbitrary functions that satisfy
the assumptions in \Cref{G105}.  Only at the very end of the proof
of (iii$'$)$\Leftrightarrow$(iii) we choose $g_1,g_2$ as in \cref{G103} for the proof
of \Cref{G14,G15}, and we use $g_1,g_2$ as in \cref{G104} for the additions to these theorems.
Again let $H_s$, $s>0$, be the trace-normalised reparameterisation of $\genResc sH$, cf.\ \Cref{G110}.
Moreover, recall that $\lambda$ denotes the Lebesgue measure.


The modified variant of \Enumref{3} reads as follows.
\begin{Ilist}
\item
	In \Cref{G14}:
	\begin{enumdash}
	\setcounter{counter_a}{2}
	\item
		For all $T\in(0,\infty)$, all $\gamma\in[0,1)$ and
		all closed, disjoint sets $A,B\subseteq\bb T$,
		the following limit relations hold:
		\begin{align}
			& \lim_{s\to 0}\lambda\bigl((0,T)\cap\sigma_{H_s}^{-1}([0,\gamma])\bigr)=0,
			\label{G44}
			\\[1ex]
			& \lim_{s\to 0}\Bigl[
			\lambda\bigl((0,T)\cap\zeta_{H_s}^{-1}(A)\bigr)
			\cdot\lambda\bigl((0,T)\cap\zeta_{H_s}^{-1}(B)\bigr)
			\Bigr]=0.
			\label{G30}
		\end{align}
	\end{enumdash}
\item
	In \Cref{G15}:
	\begin{enumdash}
	\setcounter{counter_a}{2}
	\item
		For each $T\in(0,\infty)$ there exists a sequence $(s_n)_{n\in\bb N}$ with $s_n\to 0$
		such that, for all $\gamma\in[0,1)$ and all closed, disjoint
		sets $A,B\subseteq\bb T$, the following limit relations hold:
		\begin{align}
			& \lim_{n\to\infty}\lambda\bigl((0,T)\cap\sigma_{H_{s_n}}^{-1}([0,\gamma])\bigr)=0,
			\label{G45}
			\\[1ex]
			& \lim_{n\to\infty}\Bigl[
			\lambda\bigl((0,T)\cap\zeta_{H_{s_n}}^{-1}(A)\bigr)
			\cdot\lambda\bigl((0,T)\cap\zeta_{H_{s_n}}^{-1}(B)\bigr)
			\Big]=0.
			\label{G31}
		\end{align}
	\end{enumdash}
\end{Ilist}

\noindent
Note that the statement of (iii$'$) depends on the choice of the
functions $g_1,g_2$ in \Cref{G105} because the family $(H_s)_{s>0}$
depends on $g_1$ and $g_2$.

The implication ``(ii)$\Rightarrow$(iii$'$)'' is a consequence of \Cref{G19,G23}.

\begin{proof}[Proof of \textup{(ii)}$\Rightarrow$\textup{(iii$'$)} in \Cref{G14,G15}]
	Let us first consider the situation in \Cref{G14}.
	Assume that $\lim_{t\to 0}d(H,t)=0$.  Then $\lim_{s\to0}d(H_s,T)=0$ by \cref{G113}.
	Hence, \Cref{G19,G23} applied to $H_s$ yield (iii$'$).

	Now let us consider the situation in \Cref{G15}.
	Assume that $\liminf_{t\to0}d(H,t)=0$.  By \cref{G114} there exist $s_n>0$
	with $s_n\to0$ such that $\lim_{n\to\infty}d(H_{s_n},T)=0$.
	We can apply \Cref{G19,G23} to $H_{s_n}$ to obtain (iii$'$).
\end{proof}

\noindent
The implication ``(iii$'$)$\Rightarrow$(i)'' is a consequence of \Cref{G3}.

\begin{proof}[Proof of \textup{(iii$'$)}$\Rightarrow$\textup{(i)} in \Cref{G14,G15}]
\hfill\\[0.5ex]
	By \cref{G115} and \cref{G109} we have
	\begin{equation}\label{G116}
		q_{H_s}(z) = q_{\genResc sH}(z)
		= \frac{g_3(s)}{g_2(s)}q_H\bigl(g_3(s)z\bigr).
	\end{equation}
	It follows from \Cref{G35} that
	\begin{equation}\label{G117}
		\Big|q_{H_s}\Bigl(\frac i8\Bigr)\Big|
		= \Big|q_{\genResc sH}\Bigl(\frac i8\Bigr)\Big| \asymp 1,
		\qquad s\in(0,1].
	\end{equation}
	Thus the constant Hamiltonians $\Gamma(0,1)$ and $\Gamma(0,-1)$,
	where $\Gamma$ is defined in \Cref{G32}, cannot be limit points of $(H_s)_{s\in(0,1]}$
	since $q_{\Gamma(0,1)}(z)=\infty$ and $q_{\Gamma(0,-1)}(z)=0$ by \Cref{G28}.
	Relations \cref{G116} and \cref{G117} imply that
	\begin{equation}\label{G29}
		\frac{\Im q_H\bigl(g_3(s)\frac i8\bigr)}{\big|q_H\bigl(g_3(s)\frac i8\bigr)\big|}
		= \frac{\Im q_{H_s}\bigl(\frac i8\bigr)}{\big|q_{H_s}\bigl(\frac i8\bigr)\big|}
		\asymp \Im q_{H_s}\Bigl(\frac i8\Bigr),
		\qquad s\in(0,1].
	\end{equation}
	If $T\in(0,\infty)$ and $s_n\in(0,1]$ with $s_n\to 0$ are such
	that \cref{G45} and \cref{G31} hold for all $\gamma\in[0,1)$ and
	all closed, disjoint sets $A,B\subseteq\bb T$, then
	\begin{equation}\label{G47}
		\LP_w\bigl(\rho_T(H_{s_n})\bigr)_{n\in\bb N}\subseteq\scHamInt T
	\end{equation}
	by \Cref{G3}; recall that $\rho_T:\NHam\to\NHamInt T$ is the restriction map
	\begin{Ilist}
	\item
		First assume that (iii$'$) in the sequence variant for \Cref{G15} holds.
		Then, for each $T\in(0,\infty)$, we can choose a sequence $(s_n)_{n\in\bb N}$
		that satisfies \cref{G45} and \cref{G31} for all $\gamma\in[0,1)$
		and all closed, disjoint sets $A,B\subseteq\bb T$, and hence \cref{G47}.
		\Cref{G25}\,\Enumref{2} implies that $\LP(H_s)_{s\in(0,1]}\cap\scHam\ne\emptyset$.
		Since $\Gamma(0,\pm 1)\notin\LP(H_s)_{s\in(0,1]}$,
		we find $\zeta\in\bb T\setminus\{\pm 1\}$ and some sequence $s_n\to 0$
		such that $\lim_{n\to\infty}H_{s_n}=\Gamma(1,\zeta)$.
		By the continuity of the mapping $H\mapsto q_H$ and \Cref{G28} this implies that
		\[
			\lim_{n\to\infty}q_{H_{s_n}}\Bigl(\frac i8\Bigr)
			= q_{\Gamma(1,\zeta)}\Bigl(\frac i8\Bigr)
			= \frac{\Im\zeta}{1-\Re\zeta}\in\bb R,
		\]
		and hence $\lim_{n\to\infty}\Im q_{H_{s_n}}(i/8)=0$.
		By the assumptions in \Cref{G105}, $g_3$ is continuous,
		and $g_3(s)\to\infty$ as $s\to0$.
		With $r_n\DE g_3(s_n)/8$ it follows from \cref{G29} that
		\[
			\lim_{n\to\infty}\frac{\Im q_H(ir_n)}{|q_H(ir_n)|} = 0,
		\]
		which shows \Enumref{1} in \Cref{G15}.
	\item
		Assume that (iii$'$) in the continuous variant for \Cref{G14} holds.
		We start with an arbitrary sequence $r_n\to\infty$.
		Since $g_3(s)\to\infty$ as $s\to0$ and $g_3$ is continuous, we find $s_n>0$
		for large enough $n$ such that $s_n\to0$ and $r_n=g_3(s_n)/8$.
		By \cref{G44} and \cref{G30} the relations \cref{G45} and \cref{G31}
		hold for every $T\in(0,\infty)$, every $\gamma\in[0,1)$ and
		all closed, disjoint $A,B\subseteq\bb T$ for the sequence $(s_n)_{n\in\bb N}$.
		Thus \cref{G47} holds for all $T\in(0,\infty)$, and \Cref{G25}\,\Enumref{1} gives
		$\LP(H_{s_n})_{n\in\bb N}\subseteq\scHam$.  Using that $\NHam$ is compact
		and that $\Gamma(0,\pm 1)$ cannot occur as a limit point,
		we find a subsequence $(H_{s_{n_k}})_{k\in\bb N}$ and $\zeta\in\bb T\setminus\{\pm 1\}$,
		such that $\lim_{k\to\infty}H_{s_{n_k}}=\Gamma(1,\zeta)$
		and hence $\lim_{k\to\infty}q_{H_{s_{n_l}}}(i/8)\in\bb R$ as above.
		Now relation \cref{G29} implies that
		\[
			\lim_{k\to\infty}\frac{\Im q_H(ir_{n_k})}{|q_H(ir_{n_k})|}=0.
		\]
		Since the sequence $(r_n)_{n\in\bb N}$ with $r_n\to\infty$ was arbitrary,
		the desired relation \cref{G92} follows.
	\end{Ilist}
	This finishes the proof of (iii$'$)$\Rightarrow$(i).
\end{proof}

\noindent
Showing that (iii) and (iii$'$) are equivalent is elementary.

\begin{proof}[Proof of \textup{(iii)}$\Leftrightarrow$\textup{(iii$'$)} and
of \textup{(iv)}$\Leftrightarrow$\textup{(iii$'$)} under the assumption
of \textup{\cref{G99}, \cref{G100}}]
	We prove the asserted equivalences for the continuous variant in \Cref{G14}.
	The proof for the sequence variant in \Cref{G15} is\,---\,word by word\,---\,the same.
	We proceed in several steps.  In the first two steps we show that (iii$'$)
	is equivalent to (iii$'''$) stated below.
	In the last step we prove that (iii$'''$) is equivalent to (iii)
	and\,---\,under the additional assumptions \cref{G99}, \cref{G100}\,---\,also
	equivalent to (iv).
	\begin{Elist}
	\item
		We show that (iii$'$) is equivalent to an analogous condition, say (iii$''$),
		where the limit relation \cref{G30} is required to hold for all
		open arcs $V,W\subseteq\bb T\setminus\bb R$ with
		$\ov V\cap\ov W=\emptyset$ and lengths at most $\frac\pi 2$,
		instead of all closed disjoint sets $A,B\subseteq\bb T$.

		The implication (iii$'$)$\Rightarrow$(iii$''$) is of course trivial.
		To show the converse, assume we know \cref{G44} and \cref{G30} for arcs as above.
		Let $A,B\subseteq\bb T$ be closed and disjoint.
		Then we can choose open arcs $V_1,\ldots,V_n$ and $W_1,\ldots,W_m$ of lengths
		at most $\frac\pi2$ such that
		\begin{gather}
			A\subseteq\bigcup_{i=1}^nV_i, \qquad B\subseteq\bigcup_{j=1}^mW_j,
			\label{G48}
			\\[1ex]
			\bigcup_{i=1}^n\ov{V_i}\cap\bigcup_{j=1}^m\ov{W_j}=\emptyset, \qquad
			\bb R\cap A=\bb R\cap\bigcup_{i=1}^n\ov{V_i}, \qquad
			\bb R\cap B=\bb R\cap\bigcup_{j=1}^m\ov{W_j}.
			\nonumber
		\end{gather}
		If an arc $V_i$ intersects $\bb R$, we can split it into the
		two arcs $V_i\cap\bb C^+$ and $V_i\cap\bb C^-$, and the singleton $V_i\cap\bb R$;
		here $\bb C^+$ and $\bb C^-$ are the open upper and lower half-planes respectively.
		Hence, we may assume that our arcs $V_i,W_j$
		do not intersect the real axis on the cost of adding the set $\{1,-1\}$
		to the covering, i.e.\ we can write
		\begin{equation}\label{G49}
			A \subseteq \{1,-1\}\cup\bigcup_{i=1}^nV_i,\qquad
			B \subseteq \{1,-1\}\cup\bigcup_{j=1}^mW_j
		\end{equation}
		instead of \cref{G48}.

		For any Hamiltonian we have $\zeta_H^{-1}(\{1,-1\})\subseteq\sigma_H^{-1}(\{0\})$
		by the definition of $\sigma_H$.
		Hence, \cref{G44} guarantees that
		$\lim_{s\to 0}\lambda((0,T)\cap\zeta_{H_s}^{-1}(\{1,-1\}))=0$.
		We know that, for all $i\in\{1,\ldots,n\}$ and $j\in\{1,\ldots,m\}$,
		\[
			\lim_{s\to 0}\Bigl[
			\lambda\bigl((0,T)\cap\zeta_{H_s}^{-1}(V_i)\bigr)
			\cdot\lambda\bigl((0,T)\cap\zeta_{H_s}^{-1}(W_j)\bigr)
			\Bigr]=0,
		\]
		and we obtain from \cref{G49} that also
		\[
			\lim_{s\to 0}\Bigl[
			\lambda\bigl((0,T)\cap\zeta_{H_s}^{-1}(A)\bigr)
			\cdot\lambda\bigl((0,T)\cap\zeta_{H_s}^{-1}(B)\bigr)
			\Bigr]=0.
		\]
	\item
		We make a transformation to pass from the unit circle to the real line.
		Consider the function
		\[
			\phi_+ \DF
			\left\{
			\begin{array}{rcl}
				(0,\infty) & \to & \bb T\cap\bb C^+
				\\
				x & \mapsto & \frac{1-x}{1+x}+i\sqrt{1-\bigl(\frac{1-x}{1+x}\bigr)^2}.
			\end{array}
			\right.
		\]
		This is a differentiable homeomorphism from $(0,\infty)$ onto $\bb T\cap\bb C^+$,
		and open intervals in $(0,\infty)$ correspond to open arcs in $\bb T\cap\bb C^+$.
		Moreover, for an interval $I\subseteq(0,\infty)$ we have
		\[
			\inf I=0 \;\Leftrightarrow\; 1\in\ov{\phi_+(I)}
			\qquad\text{and}\qquad
			\sup I=\infty \;\Leftrightarrow\; -1\in\ov{\phi_+(I)}.
		\]
		For any $z\in\bb T\cap\bb C^+$ the relation
		\[
			\Re\phi_+\Bigl(\frac{1-\Re z}{1+\Re z}\Bigr)=\Re z
		\]
		holds and hence also
		\[
			\phi_+\Bigl(\frac{1-\Re z}{1+\Re z}\Bigr)=z.
		\]
		Let $t\in(0,\infty)$ and assume that $\Im\zeta_H(t)>0$.
		Then $\varphi_H(t)\in(0,\frac\pi2)$, and by \cref{G94}
		we have $\pi_H(t)=\tan^2\varphi_H(t)$.
		Since $\zeta_H(t)\in\bb T\cap\bb C^+$, we have
		\[
			\zeta_H(t) = \phi_+\biggl(\frac{1-\Re\zeta_H(t)}{1+\Re\zeta_H(t)}\biggr)
			= \phi_+\biggl(\frac{1-\cos(2\varphi_H(t))}{1+\cos(2\varphi_H(t))}\biggr)
			= \phi_+\bigl(\tan^2\varphi_H(t)\bigr)
			= \phi_+\bigl(\pi_H(t)\bigr).
		\]
		Thus, for every open arc $V\subseteq\bb T\cap\bb C^+$,
		\[
			\zeta_H^{-1}(V)=\pi_H^{-1}\bigl(\phi_+^{-1}(V)\bigr).
		\]
		For the lower half-plane we proceed analogously.  Consider the function
		\[
			\phi_- \DF
			\left\{
			\begin{array}{rcl}
				(-\infty,0) & \to & \bb T\cap\bb C^-
				\\
				x & \mapsto & \frac{1-|x|}{1+|x|}-i\sqrt{1-\bigl(\frac{1-|x|}{1+|x|}\bigr)^2},
			\end{array}
			\right.
		\]
		which is a differentiable homeomorphism from $(-\infty,0)$ onto $\bb T\cap\bb C^-$
		such that open intervals in $(0,\infty)$ correspond to open arcs
		in $\bb T\cap\bb C^-$, and that, for an interval $I\subseteq(-\infty,0)$, we have
		\[
			\sup I=0 \;\Leftrightarrow\; 1\in\ov{\phi_-(I)}
			\qquad\text{and}\qquad
			\inf I=-\infty \;\Leftrightarrow\; -1\in\ov{\phi_-(I)}.
		\]
		As above one shows that, for an arc $V\subseteq\bb T\cap\bb C^-$,
		\[
			\zeta_H^{-1}(V)=\pi_H^{-1}\bigl(\phi_-^{-1}(V)\bigr).
		\]
		Now we combine the mappings $\phi_+$ and $\phi_-$;
		let $\phi\DF\bb R\setminus\{0\}\to\bb T\setminus\bb R$ be defined by
		$\phi|_{(0,\infty)}=\phi_+$ and $\phi|_{(-\infty,0)}=\phi_-$.
		The above considerations show that (iii$'$) is equivalent to
		the following condition (iii$'''$).
		\begin{enumtripledash}
		\setcounter{counter_b}{2}
		\item
			For all $T\in(0,\infty)$, all $\gamma\in[0,1)$ and
			all open intervals $I,J\subseteq\bb R\setminus\{0\}$
			with $\ov I\cap\ov J=\emptyset$ and at least one of them bounded
			the following limit relations hold:
			\begin{align}
				& \lim_{s\to 0}\lambda\bigl((0,T)\cap\sigma_{H_s}^{-1}([0,\gamma])\bigr)=0,
				\label{G118}
				\\[1ex]
				& \lim_{s\to 0}\Bigl[
				\lambda\bigl((0,T)\cap\pi_{H_s}^{-1}(I)\bigr)
				\cdot\lambda\bigl((0,T)\cap\pi_{H_s}^{-1}(J)\bigr)
				\Bigr]=0.
				\label{G119}
			\end{align}
		\end{enumtripledash}
	\item
		It follows from \cref{G51}, \cref{G111} and \cref{G108} that
		\begin{align}
			\sigma_{H_s}(t) &= \sigma_{\genResc sH}\bigl(\tau_s^{-1}(t)\bigr)
			= \sigma_H\bigl(s\tau_s^{-1}(t)\bigr),
			\label{G120}
			\\[1ex]
			\pi_{H_s}(t) &= \pi_{\genResc sH}\bigl(\tau_s^{-1}(t)\bigr)
			= \frac{g_2(s)}{g_1(s)}\pi_H\bigl(s\tau_s^{-1}(t)\bigr).
			\label{G121}
		\end{align}
		\begin{Ilist}
		\item
			To show (iii)$\Leftrightarrow$(iii$'''$), let us choose $g_1,g_2$
			as in \cref{G103}.  Then $\tau_s(t)=\mf t_s(t)$ for all $t\in(0,\infty)$,
			and \cref{G121} can be simplified to $\pi_{H_s}(t)=\pi_{H,s}(\mf t_s^{-1}(t))$.
			This and \cref{G120} show that the following equivalences hold:
			\begin{align*}
				x\in\sigma_{H_s}^{-1}([0,\gamma])
				\quad&\Leftrightarrow\quad
				x\in\mf t_s\bigl(\tfrac 1s\sigma_H^{-1}([0,\gamma])\bigr),
				\\[1ex]
				x\in\pi_{H_s}^{-1}(I)
				\quad&\Leftrightarrow\quad
				x\in\mf t_s\bigl(\pi_{H,s}(I)\bigr).
			\end{align*}
			This settles the equivalence (iii)$\Leftrightarrow$(iii$'''$).
		\item
			Finally, assume that \cref{G99} and \cref{G100} in the addition
			to \Cref{G14} hold.  Let us choose $g_1,g_2$ as in \cref{G104}.
			Then $\tau_s(t)=\frac1s\bigl(m_1(st)+m_2(st)\bigr)=t$ by \cref{G99}
			and hence $\sigma_{H_s}(t)=\sigma_H(st)$.
			For fixed $T,s\in(0,\infty)$ and $\gamma\in[0,1)$ we have
			\begin{align*}
				& \lambda\bigl((0,T)\cap\sigma_{H_s}^{-1}([0,\gamma])\bigr)
				= \lambda\bigl(\bigl\{x\in(0,T):sx\in\sigma_H^{-1}([0,\gamma])\bigr\}\bigr)
				\\[1ex]
				&= \int_{(0,T)}\mathds{1}_{\sigma_H^{-1}([0,\gamma])}(sx)\DD x
				= \frac1s\int_{(0,sT)}\mathds{1}_{\sigma_H^{-1}([0,\gamma])}(\xi)\DD\xi
				= \frac1s\lambda\bigl((0,sT)\cap\sigma_H^{-1}([0,\gamma])\bigr).
			\end{align*}
			Hence, for fixed $\gamma\in[0,1)$, the following equivalences hold:
			\begin{align*}
				\forall T\in(0,\infty)\DP \cref{G118} \ \text{is true}
				\quad&\Leftrightarrow\quad
				\forall T\in(0,\infty)\DP
				\lim_{s\to0}
				\biggl[\frac1s\lambda\Bigl((0,sT)\cap\sigma_H^{-1}([0,\gamma])\Bigr)\biggr]=0
				\\[1ex]
				&\Leftrightarrow\quad
				\lim_{t\to0}
				\biggl[\frac1t\lambda\Bigl((0,t)\cap\sigma_H^{-1}([0,\gamma])\Bigr)\biggr]=0.
			\end{align*}
			In a similar way one shows that \cref{G119} is true for every $T\in(0,\infty)$
			if and only if \cref{G66} holds.
			This establishes the equivalence of (iii$'''$) and (iv)
			and finishes the proof of \Cref{G14,G15} and their additions.
		\end{Ilist}
	\end{Elist}
	\vspace*{-4ex}
\end{proof}

\section{Hamiltonians with regularly varying diagonal}
\label{G98}

As a class of examples we consider Hamiltonians whose primitive $M$ has regularly varying
diagonal entries.
Recall that a function $f\DF(0,\infty)\to(0,\infty)$ is called
\emph{regularly varying with index $\rho$ at $0$} if
\[
	\forall t>0\DP
	\lim_{s\to 0}\frac{f(st)}{f(s)}=t^\rho;
\]
see, e.g.\ \cite[\S1.4.2]{bingham.goldie.teugels:1989}.
Typical examples of regularly varying functions are
$f(t)=t^\rho\cdot|\log t|^{\beta_1}\cdot(\log|\log t|)^{\beta_2}$
with $\rho,\beta_1,\beta_2\in\bb R$, where higher iterates of logarithms can be added.
In the theorem below we show that a Hamiltonian with regularly varying diagonal
primitives is well behaved in the sense that $d(H,t)\gtrsim 1$ unless its diagonal entries
are of the same size on the power scale, i.e.\ their indices coincide.
This is closely related to our forthcoming paper \cite{langer.pruckner.woracek:asysupp},
where we investigate Hamiltonians whose Weyl coefficients have regularly varying asymptotics
towards $+i\infty$.

\begin{theorem}\label{G67}
	Let $H$ be a Hamiltonian defined on the interval $(0,\infty)$ and assume that
	neither $h_1$ nor $h_2$ vanishes a.e.\ on some neighbourhood of the left endpoint $0$.
	Assume that $m_1$ and $m_2$ are regularly varying at $0$ with positive
	indices $\rho_1$ and $\rho_2$ respectively.  Then
	\[
		\liminf_{t\to 0}d(H,t)
		\ge 1-\biggl(\frac{\sqrt{\rho_1\rho_2}}{\frac 12(\rho_1+\rho_2)}\biggr)^2.
	\]
\end{theorem}

\begin{proof}
	Let $(\genResc sH)_{s>0}$ be the family of rescaled Hamiltonians as in \Cref{G105}
	with $g_1,g_2$ from \cref{G103},
	and let $(H_s)_{s>0}$ be the corresponding trace-normalised family
	as in \eqref{G96}.
\begin{Elist}
\item
	In the first step of the proof we show that every accumulation point of $(H_s)_{s>0}$,
	for $s\to 0$, is of a special form.
	It follows from \cref{G107} that
	\begin{equation}\label{G78}
		\biggl[\int_0^t(\genResc sH)(x)\DD x\biggr]_{ii} = \frac{m_i(st)}{m_i(s)},
		\qquad i\in\{1,2\},
	\end{equation}
	where $[C]_{ii}$ denotes the $i$th entry on the diagonal of a matrix $C$,
	and hence
	\[
		\mf t_s(t) = \tau_s(t) = \int_0^t\tr(\genResc sH)(x)\DD x
		= \frac{m_1(st)}{m_1(s)} + \frac{m_2(st)}{m_2(s)}
	\]
	where $\mf t_s$ and $\tau_s$ are defined in \cref{G76} and \cref{G111} respectively.
	Set $\mf t(t)\DE t^{\rho_1}+t^{\rho_2}$ for $t\in(0,\infty)$.
	The assumptions about $m_1$ and $m_2$ and the Uniform Convergence Theorem
	for regularly varying functions
	(see, e.g.\ \cite[Theorem~1.5.2]{bingham.goldie.teugels:1989}) imply that
	$\lim_{s\to 0}\mf t_s(t)=\mf t(t)$ locally uniformly for $t\in(0,\infty)$.
	The functions $\mf t_s$ and $\mf t$ are continuous and increasing bijections
	from $(0,\infty)$ onto itself, and it follows
	that also $\lim_{s\to 0}\mf t_s^{-1}(T)=\mf t^{-1}(T)$ for all $T\in(0,\infty)$.

	Let $s_n\to 0$ be a sequence such that the limit $\wt H\DE\lim_{n\to\infty}H_{s_n}$
	exists, and let $\widehat H$ be the reparameterisation
	defined by $\widehat H\DE(\wt H\circ\mf t)\cdot \mf t'$.
	Using \cref{G78} we find, for $T\in(0,\infty)$ and $i\in\{1,2\}$, that
	\begin{align*}
		\biggl[\int_0^{\mf t^{-1}(T)}\widehat H(t)\DD t\biggr]_{ii}
		&= \biggl[\int_0^{\mf t^{-1}(T)}\wt H\bigl(\mf t(t)\bigr)\mf t'(t)\DD t\biggr]_{ii}
		= \biggl[\int_0^T\wt H(x)\DD x\biggr]_{ii}
		\\[1.5ex]
		&= \lim_{n\to\infty}\biggl[\int_0^T H_{s_n}(x)\DD x\biggr]_{ii}
		= \lim_{n\to\infty}\biggl[\int_0^{\mf t_{s_n}^{-1}(T)}(\genResc{s_n}H)(t)
		\DD t\biggr]_{ii}
		\\[1.5ex]
		&= \lim_{n\to\infty}\frac{m_1(s_n\mf t_{s_n}^{-1}(T))}{m_1(s_n)}
		= \mf t^{-1}(T)^{\rho_1},
	\end{align*}
	again by the Uniform Convergence Theorem.  Hence $\widehat H$ is of the form
	\begin{equation}\label{G75}
		\widehat H(t) =
		\begin{pmatrix}
			\rho_1 t^{\rho_1-1} & \text{\ding{100}}
			\\[1ex]
			\text{\ding{100}} & \rho_2 t^{\rho_2-1}
		\end{pmatrix}
	\end{equation}
	where the off-diagonal entries are unknown.
\item
	For Hamiltonians $\widehat H$ of the form \cref{G75} an estimate
	for $d(\widehat H,t)$ holds.
	With $\hat h_j$ being the entries of $\widehat H$ we have
	\[
		|\hat h_3(t)| \le \sqrt{\hat h_1(t)\hat h_2(t)}
		=\sqrt{\rho_1\rho_2}\,t^{\frac 12(\rho_1+\rho_2)-1}
	\]
	and hence
	\[
		|\widehat m_3(t)| \le \int_0^t |\hat h_3(x)\DD x|
		\le \frac{\sqrt{\rho_1\rho_2}}{\frac 12(\rho_1+\rho_2)}t^{\frac 12(\rho_1+\rho_2)},
	\]
	from which we find that, for all $t>0$,
	\begin{equation}\label{G81}
		d(\widehat H,t)
		\ge 1-\biggl(\frac{\sqrt{\rho_1\rho_2}}{\frac 12(\rho_1+\rho_2)}\biggr)^2.
	\end{equation}
\item
	We make a limiting argument to complete the proof.
	Let $(t_n)_{n=1}^\infty$ be a sequence of positive numbers with $t_n\to 0$.
	Fix $T>0$ and let again $u(s)$ be the function in \cref{G80}.
	For large enough $n$, choose $s_n\to 0$ such that $u(s_n)=t_n$,
	and extract a subsequence $(s_{n(k)})_{k\in\bb N}$ such that the
	limit $\widehat H\DE\lim_{k\to\infty}H_{s_{n(k)}}$ exists.
	Using \cref{G108}, \cref{G82} and \cref{G81} we obtain
	\begin{align*}
		d(H,t_{n(k)}) &= d\bigl(H,u(s_{n(k)})\bigr)
		= d\bigl(\genResc{s_{n(k)}}H,\mf t_{s_{n(k)}}^{-1}(T)\bigr)
		\\[1ex]
		&= d(H_{s_{n(k)}},T)\stackrel{k\to\infty}\longrightarrow d(\widehat H,T)
		\ge 1-\biggl(\frac{\sqrt{\rho_1\rho_2}}{\frac 12(\rho_1+\rho_2)}\biggr)^2.
	\end{align*}
	Since the $(t_n)$ was arbitrary, the claim follows.
\end{Elist}
\vspace*{-3ex}
\end{proof}

\noindent
As a consequence, if $\rho_1\ne\rho_2$ in \Cref{G67}, then \Enumref{2} in \Cref{G15}
is not satisfied and hence neither is \Enumref{1}
(under the assumption that \cref{G1} holds), i.e.\
one has $\liminf_{y\to\infty}\frac{\Im q_H(iy)}{|q_H(iy)|}>0$.
If, on the other hand, the diagonal entries themselves (and not just their primitives)
are regularly varying with the same index, then the situation is different.

\begin{proposition}\label{G102}
	Assume that $h_1(t),h_2(t)>0$ a.e., that $h_1,h_2$ are regularly varying
	with the same index $\alpha>-1$, and set $h_3(t)\DE\sqrt{h_1(t)h_2(t)}$, $t\in(0,\infty)$.
	Then $\lim_{t\to 0}d(H,t)=0$ and hence $\lim_{y\to\infty}\frac{\Im q_H(iy)}{|q_H(iy)|}=0$.
\end{proposition}
\begin{proof}
	The off-diagonal entry $h_3$ is also regularly varying with index $\alpha$.
	It follows from Karamata's Theorem
	(e.g.\ \cite[Theorem~1.5.10]{bingham.goldie.teugels:1989} transformed
	from the asymptotics at $\infty$ to the asymptotics at $0$ by a change of variable)
	that $m_i(t)=\frac{1}{1+\alpha}th_i(t)(1+\Smallo(1))$ as $t\to0$ for $i=1,2,3$.
	Hence
	\begin{align*}
		d(H,t) &= \frac{m_1(t)m_2(t)-m_3(t)^2}{m_1(t)m_2(t)}
		= \frac{h_1(t)h_2(t)\bigl(1+\Smallo(1)\bigr)-h_3(t)^2\bigl(1+\Smallo(1)\bigr)}{h_1(t)h_2(t)\bigl(1+\Smallo(1)\bigr)}
		\\[1ex]
		&= \frac{h_1(t)h_2(t)\Smallo(1)}{h_1(t)h_2(t)\bigl(1+\Smallo(1)\bigr)}
		\to 0
	\end{align*}
	as $t\to0$.  The last statement follows from \Cref{G14}.
\end{proof}

\bigskip

\noindent
\textbf{Acknowledgements.} \\
We thank a referee for suggestions that enabled us to shorten some of the proofs.



{\footnotesize
\begin{flushleft}
	M.~Langer \\
	Department of Mathematics and Statistics \\
	University of Strathclyde \\
	26 Richmond Street \\
	Glasgow G1 1XH \\
	UNITED KINGDOM \\
	email: \texttt{m.langer@strath.ac.uk} \\[5mm]
\end{flushleft}
\begin{flushleft}
	R.~Pruckner \\
	Institute for Analysis and Scientific Computing \\
	Vienna University of Technology \\
	Wiedner Hauptstra{\ss}e 8--10/101 \\
	1040 Wien \\
	AUSTRIA \\
	email: \texttt{raphael.pruckner@tuwien.ac.at} \\[5mm]
\end{flushleft}
\begin{flushleft}
	H.\,Woracek\\
	Institute for Analysis and Scientific Computing\\
	Vienna University of Technology\\
	Wiedner Hauptstra{\ss}e\ 8--10/101\\
	1040 Wien\\
	AUSTRIA\\
	email: \texttt{harald.woracek@tuwien.ac.at}
\end{flushleft}
}


%
%

%
%
%
\end{document}